\documentclass[11pt]{amsart}

\usepackage{amssymb,amsthm,amsmath}
\usepackage[numbers,sort&compress]{natbib}
\usepackage{color}
\usepackage{graphicx}
\usepackage{tikz}
\usepackage[colorlinks,
            linkcolor=blue,
            anchorcolor=blue,
            citecolor=blue
            ]{hyperref}

\newcommand{\SL}{\mathrm{SL}(2,\mathbb{R})}

\newcommand{\la}{\langle }
\newcommand{\ra}{\rangle }

\let\newpf\proof \let\proof\relax 
\newenvironment{pf}{\newpf[\proofname]}{\qed\endtrivlist}

\newcommand{\bQ}{\overline{Q}}
\newcommand{\ba}{\overline{A}}

\newcommand{\cF}{\mathcal{F}}

\newcommand{\CD}{\rm CD}

\def\be{\begin{equation}}
\def\ee{\end{equation}}

\def\ba{{\begin{align}}}
\def\ea{{\end{align}}}

\def\bm{\begin{matrix}}
\def\em{\end{matrix}}

\def\SL{{\mathrm{SL}}}

\def\0{{\mathbf 0}}

\def\cal{\mathcal}

\newtheorem{Theorem}{Theorem}[section]
\newtheorem{Lemma}{Lemma}[section]
\newtheorem{Proposition}{Proposition}[section]
\newtheorem{Corollary}{Corollary}[section]
\newtheorem{Remark}{Remark}[section]

\newtheorem{Definition}{Definition}[section]

\numberwithin{equation}{section}

\def \bn {\hfill \\ \smallskip\noindent}

\theoremstyle{definition}

\def\proof{\bn {\bf Proof.} }

\newcommand{\C}{{\mathbb C}}

\newcommand{\N}{{\mathbb N}}
\newcommand{\Q}{{\mathbb Q}}
\newcommand{\R}{{\mathbb R}}
\newcommand{\T}{{\mathbb T}}

\newcommand{\Z}{{\mathbb Z}}

\def\B0{{\bold{0}}}

\catcode`\@=12

\def\Empty{}
\newcommand\oplabel[1]{
  \def\OpArg{#1} \ifx \OpArg\Empty {} \else
    \label{#1}
  \fi}

\newcommand{\comm}[1]{}
\newcommand{\comment}[1]{}

\begin{document}

\title[high-dimensional base reducibility]{Absolutely Continuous Spectrum of Multifrequency Quasiperiodic Schr\"odinger operator}

\author {Xuanji Hou}
\address{
School of Mathematic and Statics, Central China Normal University, Wuhan
430079,  China } \email{hxj@mail.ccnu.edu.cn}

\author {Jing Wang}
\address{
Department of  Mathematics, School of Science, Nanjing University of Science and Technology, Nanjing 210094, China}
\email{jingwang018@gmail.com, jing.wang@njust.edu.cn}

\author{Qi Zhou}
\address{
Chern Institute of Mathematics and LPMC, Nankai University, Tianjin 300071, China}
 \email{qizhou628@gmail.com, qizhou@nankai.edu.cn}
\date{\today}

\begin{abstract}
In this paper, we prove that for any $d$-frequency analytic quasiperiodic Schr\"odinger operator, if the  frequency is weak Liouvillean, and the potential is small enough, then the corresponding operator has absolutely continuous spectrum.  Moreover, in  the case   $d=2$, we even establish the existence of ac spectrum under small
potential and some super-Liouvillean frequency,    and this result is optimal due to a recent counterexample of
Avila and Jitomirskaya \cite{AJ2}.
\end{abstract}

\date{\today}

\setcounter{tocdepth}{1}

\maketitle

\section{Introduction and main results}
This work is concerned with  one-dimensional analytic
Schr\"{o}dinger operators  $H_{v,\alpha, \phi}$ defined on $l^2(\Z)$
\begin{eqnarray}
(H_{v,\alpha, \phi}u)_n= u_{n-1} +u_{n+1} +  v (\phi+n\alpha) u_n,
\end{eqnarray}
where $v\in C^\omega(\R^d/ \Z^d, \R)$,  $\phi\in \R^d/ \Z^d=:\T^d$ and
 $\alpha \in \R^d $ is \textit{rational independent}, i.e., $\langle \mathbf{k}, \alpha\rangle \notin \Z$ for all non-zero $\mathbf{k}\in\Z^d$.
 We recall that $ v$,  $\phi$,  and $\alpha$ are the \textit{potential}, the \textit{phase} and the \textit{frequency} respectively. Because of its strong connection with physics, analytic quasiperiodic  Schr\"{o}dinger operators  have been fruitfully studied for a long time. We invite the reader to consult the recent surveys \cite{D,JM}  and the references therein.

 The spectrum set of  $H_{ v,\alpha, \phi}$ is a compact set of
 $\R$, independent of the phase $\phi$, since the base dynamics $\phi\mapsto \phi+\alpha$ is minimal, and we will denote it by  $\Sigma(v, \alpha)$ for short.
  The spectral measure of  $H_{  v,\alpha, \phi}$ can be decomposed into atomic, singular continuous and absolutely continuous (ac) parts, while ac part corresponds to the strongest transport properties.   The ac spectrum part of  $H_{ v,\alpha, \phi}$
  is also independent of the phase $\phi$ \cite{LS}, and we use   $\Sigma_{ac}(v, \alpha)$ to denote it.  We are most interested in when $ \Sigma_{ac}( v, \alpha) \neq \varnothing$. The well-known  \textit{Kotani's theory}  \cite{K, S} tells us  that $\Sigma_{ac}( v, \alpha)$ has a close relationship with the \textit{Lyapunov exponent}
\begin{eqnarray}
L(\alpha, E-  v )= \lim_{n\rightarrow \infty} \frac{1}{n} \int_{\R^d/ \Z^d} \ln \|S_{E- v}^{(n)}(\phi) \|d\phi \geq 0,
\end{eqnarray}
where $E\in \R$ is the energy and
\begin{eqnarray*}
S_{E- v}^{(n)}(\phi)= \left(
\begin{array}{ccc}
  E - v(\phi+(n-1)\alpha) &   -1\cr
   1 & 0\end{array}
   \right)\cdots \left(
\begin{array}{ccc}
 E-  v(\phi) &   -1\cr
   1 & 0\end{array}
   \right).
\end{eqnarray*}
 Let $\Sigma_0( v, \alpha)=\{ E\in \R |  L(\alpha, E-v)=0  \}$. Then $ \Sigma_{ac}(  v, \alpha)$ is the essential closure\footnote{Given a set $A\subseteq \R$, the essential closure of $A$ is defined as the set $$\{E\in\R |\forall \epsilon>0, A\cap(E-\epsilon,E+\epsilon)\, \hbox{is of  positive Lebesgue measure} \}.$$} of  $ \Sigma_0( v, \alpha)$ \cite{K, S}.

In order to study the spectral property of $H_{v,\alpha,\phi}$, it is very useful to introduce the \textit{coupling constant} $\lambda$ in front of the potential $v$. By the symmetry, we just assume $\lambda>0$. If the coupling constant $\lambda$ is sufficiently large, then  $L(\alpha, E-\lambda v)>0$ for any $E\in \R$\cite{B2,  BG,GS, H, SS}. Therefore by Kotani's Theory, there is no hope to find $\Sigma_{ac}(\lambda v, \alpha)$ of $H_{\lambda v,\alpha, \theta}$ for large coupling constant.

 If the coupling constant $\lambda$ is sufficiently  small, there is a tendency of occurrence of  ac spectrum,  where the  arithmetic properties of the frequency come into play.   We recall that  $\alpha$ is Diopantine (denote it by $\alpha\in DC$), if there exist  $\gamma,\tau>0$, such that for any $\mathbf k\in\Z^d\backslash\{0\}$
$$\|\langle \mathbf{k}, \alpha\rangle\|_{\R/\Z} =\min_{j\in \Z}|\langle \mathbf k, \alpha\rangle-j|\geq \frac{\gamma}{ |\mathbf k|^{\tau}}.$$ Based on KAM method,  Dinaburg-Sinai \cite{DS} proved that $\Sigma_{ac}(\lambda v, \alpha) \neq \varnothing$ in the \textit{perturbative small regime $\lambda<\lambda_0$}.
Here perturbative means that $\lambda_0$ depends on $\alpha$ through the Diophantine constants $\gamma,\tau$. Also under the same assumption,  Eliasson \cite{E92} showed that in fact we have  $\Sigma(\lambda v, \alpha)=\Sigma_{ac}(\lambda v, \alpha)$.

Moreover,  in the case of $d=1$, one can even anticipate some nonperturbative results.  Making use of the specificity of one frequency,  some new elaborate   techniques are developed to prove some sharp results.  As $\alpha\in DC$, based on non-perturbative Anderson localization result, Avila-Jitomirskaya \cite{AJ1} proved that there exists $\lambda_1$ which  does not depend on $\alpha$, such that  $\Sigma(\lambda v, \alpha)=\Sigma_{ac}(\lambda v, \alpha)$ as $\lambda< \lambda_1$. Recently,  Avila-Fayad-Krikorian \cite{AFK} and Hou-You \cite{HY} independently  developed  non-standard KAM
techniques, showed that  $\Sigma_{ac}(\lambda v, \alpha) \neq \varnothing$ for  $\lambda< \lambda_1(v)$ and for any  irrational $\alpha$. Indeed, Avila \cite{Aac1,A2} was able to show that  $\Sigma(\lambda v, \alpha)=\Sigma_{ac}(\lambda v, \alpha)$.

If the coupling constant  $\lambda$ is neither too large nor too small,  which is the so-called \textit{global case}, the problem is extremely complicated, since few methods can be successfully applied in obtaining detailed informations of the spectrum.
Nevertheless,  in  the one frequency case, we would like to highlight Avila's fascinating global theory of  analytic Schr\"odinger operators \cite{A1} and his well-known spectral dichotomy conjecture(now proved by himself in \cite{Aac1,A2}).  He shows that for typical analytic one-frequency  Schr\"odinger operator, it doesn't have singular continuous spectrum. This is an already very deep result from the ac spectrum side, as was former pointed out by Avila and Krikorian \cite{AK}, $\Sigma_0(\lambda v, \alpha)$ may support non-ac spectrum. Avila's global theory \cite{Aac1, A2} thus establishes the deep relations between the existence of ac spectrum  and the vanishing of the Lyapunov exponent in a band \cite{A1}.

\subsection{Multifrequency Schr\"odinger operator}
However,  not to mention  Avila's global theory can't be generalized to the multifrequency case \cite{DK},  one can even not anticipate non-perturbative results for the multifrequency case (a counter example is given by Bourgain \cite{B1}).  Comparing to the one frequency case,  as pointed by  Eliasson \cite{Eliconm}: \textquotedblleft while almost everything is known for one frequency case, almost everything is unknown for multifrequency case\textquotedblright.

In this paper, we establish two  results which confirm the existence of  ac spectrum for multifrequency quasiperiodic  Schr\"odinger operator beyond the Diophantine frequency. Before introducing them,
we need to give some necessary definitions.

For $\tilde \alpha\in \R\backslash \Q$, $\alpha' \in\T^{d-1}$,  we  say that the pair $(\tilde\alpha, \alpha')$ is  \textit{weak Liouvillean}, if there exist $\gamma, \tau>0, 0< \tilde{U}<\infty$, such that
\[\widetilde{U}=\widetilde{U}(\tilde\alpha)=:\sup_{n>0} \frac{\ln\ln \tilde q_{n+1}}{\ln \tilde q_n}, \]
\begin{eqnarray*}
\label{rev200329-wldef}\| k \tilde \alpha+ \la l,\alpha'\ra\|_{\R/\Z} \geq \frac{\gamma}{(|k|+|l|)^\tau}\
\ \textrm{for} \  k\in \Z, l\in \Z^{d-1}\backslash\{0\},
\end{eqnarray*}
where $\frac{\tilde p_n}{\tilde q_n}$  is  the continued fraction approximates to $\tilde\alpha$. 
We will  denote by $WL(\gamma, \tau, \tilde{U})$ the set of such pairs and by $WL$ the union
\[WL=\bigcup_{\gamma, \tau>0, 0< \tilde{U}<\infty}WL(\gamma, \tau, \tilde{U}). \]
 It is obvious that  $WL$ is of full Lebesgue
measure. Here we remark that for $(\tilde\alpha,\alpha')\in DC$, the growth of $\{\tilde q_n\}_{n\in\N}$ is no more than polynomial, while for $(\tilde \alpha,\alpha')\in WL$, the growth of  $\{\tilde q_n\}_{n\in\N}$ can be super-exponential.  One of our main results is the following:

\begin{Theorem}\label{thm-ana}
Assume that $\alpha= (\tilde\alpha, \alpha')\in WL $, $v\in C^\omega(\T^d,\R)$. Then there exists $\lambda_2=\lambda_2(\alpha,v)>0$ such that if $\lambda<\lambda_2$, then $\Sigma_{ac}(\lambda v, \alpha)$ is of positive Lebesgue measure.
\end{Theorem}

For any rational independent $\alpha\in \T^d$, in order to measure how Liouvillean $\alpha$ is,  one very useful definition is the following:
\begin{equation}\label{beta}
\beta(\alpha):= \limsup_{\mathbf k\in \Z^{d}} \frac{1}{ |\mathbf k|} \ln \frac{1}{\|\la \mathbf k,\alpha\ra \|_{\R/\Z}}.
\end{equation}
One can show that $WL\cap\{\alpha\in\T^d\,\big|\, \beta(\alpha)=+\infty\}=\varnothing$, and for any $\beta_*\in [0,+\infty)$ we have $WL\cap\{\alpha\in\T^d\,\big|\, \beta(\alpha)=\beta_*\}\neq \varnothing$.
Then as a direct corollary of Theorem \ref{thm-ana}, we have the following:

\begin{Corollary}\label{thm-ana-plus}
For any $0\leq \beta_*<\infty$, there exists $\alpha \in \T^d$ with $\beta(\alpha)=\beta_*$, and for any  $v\in C^\omega(\T^d,\R)$, there exists $\lambda_3=\lambda_3(\alpha,v)>0$, such that  if $\lambda<\lambda_3$, then  $\Sigma_{ac}(\lambda v, \alpha)$ is of positive Lebesgue measure.
\end{Corollary}

\subsection{Two-dimensional frequency Schr\"odinger operator}

By Corollary \ref{thm-ana-plus}, we establish the existence of ac spectrum for non-super-Liouvillean frequency, i.e., $\beta(\alpha)<\infty$, yet if the frequency is two dimensional, we are able to show more:

\begin{Theorem}\label{thm-2-fre}
There exists $\alpha\in \T^2$ satisfying $\beta(\alpha)=\infty$,  and for any  $v\in C^\omega(\T^2,\R)$, there exists $\lambda_4=\lambda_4(\alpha,v)>0$,  such that if $\lambda<\lambda_4$, then $\Sigma_{ac}(\lambda v, \alpha)$ is of positive Lebesgue measure.
\end{Theorem}

If the frequency is two dimensional, we can only prove the existence of ac spectrum for \textit{some} super-Liouvillean frequency (i.e. $\beta(\alpha)=\infty$). One thus wonder whether one can prove this for \textit{all} super-Liouvillean frequency, as in the one-frequency  case \cite{Aac1,A2}. However, this is  impossible due to
recent remarkable counterexample of  Avila-Jitomirskaya \cite{AJ2}, who
proved that there exists $\alpha\in \T^2$ satisfying $\beta(\alpha)=\infty$, such that for  typical  analytic potential $v\in C^\omega(\T^2,\R)$, $\Sigma_{ac}(v, \alpha)$ is empty.  One may understand this result partially in the following way: in establishing the existence of ac spectrum for any irrational frequency, one key ingredient is the Denjoy-Koksma lemma \cite{AFK}; however, this estimate fails in the multifrequency case, as pointed out by Yoccoz \cite{Y}.

We also want to point out that, our result (Theorem \ref{thm-2-fre}) and Avila-Jitomirskaya's \cite{AJ2} are both constructive. In section \ref{sec-2-freq-construction}, we will show a more precise version of this result, and construct the frequency very explicitly.    Readers are invited to consult Remark \ref{compareaj} for more discussions, and we will show that the existence of ac spectrum  depends on quite dedicate arithmetic properties of the frequency.

\section{Preliminaries}

In this section, some frequently used  facts and notations are introduced.

\subsection{Notations on functions}
An integrable  real or complex valued function $f$ on $\T^d$ has the Fourier expansion $f=\sum_{\mathbf k\in\Z^d}\hat{f}(\mathbf k)e^{2\pi i
\langle \mathbf k, \phi\rangle }$ with $\hat{f}(\mathbf k)=\int_{\T^d}f(\phi)e^{- 2\pi i
\langle \mathbf k, \phi\rangle }d\phi$. For any $N>0$, $\mathcal{T}_N$ and $\mathcal{R}_N$  are used
to denote the truncation operators:
\begin{equation}
\mathcal{T}_N (f)=\sum_{|\mathbf k|<N} \hat f(\mathbf k)e^{2\pi i\la \mathbf k,\phi\ra},
\  \mathcal{R}_N (f) =\sum_{ |\mathbf k|\geq N} \hat f(\mathbf k) e^{2\pi i
\la \mathbf k,\phi\ra}.
\end{equation}
Let $\mathfrak r>0$, and we denote by $C^\omega_{\mathfrak r}(\T^d,\R)$ ($C^\omega_{\mathfrak r}(\T^d,\C)$) the
set of all  real (complex) valued functions admitting an analytic extension on
\begin{equation}
\T^d_{\mathfrak r}:= \{ \phi\in \T^d\  \big | \  |\Im \phi_1|\leq \mathfrak r,\cdots, |\Im \phi_d|\leq  \mathfrak r \}.
\end{equation}
For any $f\in C^\omega_{\mathfrak r}(\T^d,\R)$ ($f\in C^\omega_{\mathfrak r}(\T^d,\C)$), there is
\begin{equation}
\|f\|_{\mathfrak r} :=  \sup_{ \phi\in \T^d_{\mathfrak r} } | f(\phi)|<\infty,
\end{equation}
and it is well-known that $|\hat{f}(\mathbf k)|\leq \|f\|_{\mathfrak r}e^{-2\pi |\mathbf k|\mathfrak r}$ for all $\mathbf k\in\Z^d$.

In this article, we also  frequently consider real (complex) valued functions admitting an analytic extension on
\begin{equation}
\T^d_{r,s}:= \{ \phi\in \T^d\  \big | \  |\Im \phi_1|\leq r,|\Im \phi_2|\leq s,\cdots, |\Im \phi_d|\leq  s \},
\end{equation}
where $r,s>0$, and use $C^\omega_{r,s}(\T^d,\R)$ ($C^\omega_{r,s}(\T^d,\C)$) to denote the set of all such functions.
Any $f$ in $C^\omega_{r,s}(\T^d,\R)$ ($C^\omega_{r,s}(\T^d,\C)$) satisfies
\begin{equation}
\|f\|_{r,s} :=  \sup_{ \phi\in \T^d_{r,s} } | f(\phi)|<\infty,
\end{equation}
and $|\hat{f}(k,l)|\leq \|f\|_{r,s}e^{-2\pi (|k|r+|l|s)}$ for all $\mathbf (k,l)\in\Z\times \Z^{d-1}$.

Let $M_2$ denote the set of all $2$ by $2$ complex matrices and $\|\cdot\|$ denote the matrix norm.
An integrable $M_2$ valued function $F$ on $\T^d$ has the Fourier expansion $F=\sum_{\mathbf k\in\Z^d}\hat{F}(\mathbf k)e^{2\pi i
\langle \mathbf k, \phi\rangle }$ with $\hat{F}(\mathbf k)=\int_{\T^d}F(\phi)e^{- 2\pi i
\langle \mathbf k, \phi\rangle }d\phi$. For any $N>0$, we define truncation operators $\mathcal{T}_N$ and $\mathcal{R}_N$
in similar way.  Let $C^\omega_{\mathfrak r}(\T^d,\, M_2)$ denotes  the
set of all  $M_2$-valued functions on $\T^d$ admitting an analytic extension on $\T^d_{\mathfrak r}$. For any $F\in C^\omega_{\mathfrak r}(\T^d,\, M_2)$  we define
\begin{equation}
\|F\|_{\mathfrak r} :=  \sup_{ \phi\in \T^d_{\mathfrak r} } \| F(\phi)\|,
\end{equation}
 and there are estimates $\|\hat{F}(\mathbf k)\|\leq \|F\|_{\mathfrak r}e^{-2\pi |\mathbf k|\mathfrak r}\, (\mathbf k\in\Z^d)$.
 We also use $C^\omega_{\mathfrak r}(\T^d,\, *)$ to denote the set of $*$-valued functions in $C^\omega_{\mathfrak r}(\T^d,\, M_2)$,
 where $*$ can be $sl(2,\R)$,
$SL(2,\R)$,  $PSL(2,\R)$\footnote{$PSL(2,\R)$ denotes the quotient group $SL(2,\R)/\{\pm I\}$.}  etc..  Similarly, for $r,s>0$, one can  also define
$C^\omega_{r,s}(\T^d,\, *)$ with corresponding norm $\|\cdot\|_{r,s}$ on it.

%

\subsection{Continued fraction expansion and CD bridge}

Let $\alpha\in\R^1$ be irrational. Define $a_0=[\alpha],
\alpha_{0}=\alpha-a_0,$ and inductively for $j\geq 1$,
$$a_j=[\alpha_{j-1}^{-1}],\qquad \alpha_j=\{ \alpha_{j-1}^{-1}\}=\alpha_{j-1}^{-1}-a_j.$$
We define $p_0=0,  p_1=1, q_0=1, q_1=a_1$ and inductively,
\begin{eqnarray*}
 p_j=a_jp_{j-1}+p_{j-2}, \quad  q_j=a_jq_{j-1}+q_{j-2}.
\end{eqnarray*}
Then  the sequence $(q_n)$ satisfies
\begin{eqnarray*}
&& \forall 1 \leq q < q_n,\quad \|q\alpha\|_{\R/\Z} \geq
\|q_{n-1}\alpha\|_{\R/\Z},\\
&& \frac{1}{q_n+q_{n+1}}\leq \|q_n \alpha\|_{\R/\Z}\leq {1 \over
q_{n+1}}.
\end{eqnarray*}$\,$

Let us now introduce the CD bridge.

\begin{Definition}[{\cite{AFK}}]\label{CDbridge}
Let $\mathcal A, \mathcal B, \mathcal C\in\R$ with $0<\mathcal A \leq \mathcal B\leq \mathcal C$. We say that the pair of denominators
$(q_l,q_n)$ forms a $\CD(\cal A,\cal B,\cal C)$ bridge if
\begin{itemize}
\item $q_{j+1}\leq q_j^{ \cal A}, \quad \forall j=l,\ldots,n-1$
\item $q_l^{\cal C}\geq q_n\geq q_l^{\cal B}$.
\end{itemize}
\end{Definition}

 In the following, for simplicity, we will fix a subsequence $(q_{n_j})$ of
$(q_n)$, denoted by $(Q_j)$, and the  subsequence $(q_{n_j+1})$,
denoted by $(\bQ_j)$.

\begin{Lemma}[{\cite{AFK}}]\label{CDbridge}
For any ${\cal A}>0$, there exists  a subsequence $(Q_j)$ such that $Q_0=1$ and for each $j\geq 0$,
$Q_{j+1}\leq \bQ_j^{{\cal A}^4}$. Furthermore, either $\bQ_j\geq
Q_j^{\cal A}$, or the pairs $(\bQ_{j-1},Q_{j})$ and $(Q_j,Q_{j+1})$
are both $\CD({\cal A},{\cal A},{\cal A}^3)$ bridges.
\end{Lemma}

\begin{Corollary}[{\cite{KWYZ}}]\label{cd-coroll}
Let $\alpha\in \R\backslash\Q$. If $\widetilde{U}(\alpha)<\infty,$ then for the subsequence $(Q_j)$ selected in Lemma \ref{CDbridge} we have $Q_j\geq
Q_{j-1}^\cal A$ for every $j\geq 1$. Furthermore, we have
$$\sup_{j>0} \frac{\ln\ln Q_{j+1}}{\ln Q_j}\leq U(\alpha),$$ where
$U(\alpha):=\widetilde{U}(\alpha)+4\frac{\ln \mathcal A}{\ln 2}<\infty$.
\end{Corollary}

 \subsection{Cocycles, Fibered rotation number}\label{sec-cocycles}

Let $\alpha\in \T^d$ be rational independent and $G=SL(2,\R),SO(2,\R), etc.$. A quasi-periodic $G$ cocycle is a pair $(\alpha, A)\in\T^d\times C^r(\T^d, G)$ ($r=0,1,\cdots,\infty,\omega$), which represents
 \begin{eqnarray*}
\quad (\alpha, \, A): \, \T^d\times\R^2\rightarrow\T^d\times\R^2,\, (\phi,x)\mapsto(\phi+\alpha,\,A(\phi)x).
\end{eqnarray*}
We consider the iterations of  $(\alpha, A)$ which can be written as
\[(\alpha, A)^{\circ n}=\underbrace{(\alpha, A)\circ \cdots \circ(\alpha, A)}_{\hbox{$n$ times}}=(n\alpha, A^{(n)})\]
where
\[A^{(n)}(\cdot)=\left\{
\begin{array}{lll}
& A(\cdot+(n-1)\alpha)\cdots A(\cdot),\quad &\hbox{as $n>0$,}\\
 & I, \quad &\hbox{as $n=0$,}\\
 & A(\cdot+n\alpha)^{-1}\cdots A(\cdot-\alpha)^{-1},\quad &\hbox{as $n<0$}.
\end{array} \right.\]
We say that two cocycles $(\alpha, A), (\alpha, \widetilde{A})$ are   $C^r$  \textit{conjugate}, if there exists $B\in C^r(\T^d,\, PSL(2, \mathbb{R}))$, such that
$$\widetilde{A}(\cdot)=B(\cdot+\alpha)A(\cdot)
B^{-1}(\cdot).$$
Furthermore, if $\widetilde{A} \in C^r(\T^d,\, SO(2, \mathbb{R}))$
we say that $(\alpha, A)$ is  $C^r$ \textit{rotations reducible}.   It is obvious that rotations reducible cocycles are always bounded.

For any  $C^0$ quasi-periodic $SL(2,\R)$ cocycle
$(\alpha, A)$ with $A$ being homotopic to identity, it induces a projective cocycle on $\T^d\times \mathbb S^1$
 \[f_{(\alpha, A)}: \, \T^d\times\mathbb S^1\rightarrow\T^d\times\mathbb S^1,\, (\phi,\overrightarrow{e}(x))\mapsto(\phi+\alpha,\,\frac{A(\phi)\overrightarrow{e}(x)}{\|A(\phi)\overrightarrow{e}(x)\|}),
\]
with $\overrightarrow{e}(x)=(\cos2\pi x, \sin 2\pi x)^T$, which is also homotopic to identity.
Then $f_{(\alpha, A)}$ has a lift $F_{(\alpha, A)} :\ \T^d\times\R\circlearrowleft$ of the form
$(\phi, x)\mapsto(\phi+\alpha,\,\tilde f_{A}(\phi,x))$
with $\tilde f_{A}(\phi,x)=x+d_A(\phi,x)$, where $d_A$ is $\Z-$periodic in both variables, such that
$\overrightarrow{e}\left(\tilde f_A(\phi,x)\right)=\frac{A(\phi)\overrightarrow{e}(x)}{\|A(\phi)\overrightarrow{e}(x)\|}$.
As $\alpha$ has rationally independent coordinates, then the limit
\[\tilde\rho(\alpha,A):=\lim_{n\rightarrow\infty}\frac{\tilde f_A(F_{(\alpha,A)}^{n-1}(\phi,x))-x}{n}\]
exists for all $(\phi,x)\in\T^d\times\R$ and it is independent of $(\phi,x)$\cite{H}.
We call $\rho(\alpha, A):=\tilde\rho(\alpha,A)$ (mod 1) the \textit{fibered rotation number} of $(\alpha,A)$.

One can establish fibered rotation numbers for Schr\"{o}dinger  cocyles $(\alpha, S_{E-v})$, since $S_{E-v}$ is obviously homotopic to identity, where
$E\in\R$ and $S_{E-v}(\cdot)=\left(\begin{array}{cc} E-v(\cdot)& -1\\ 1& 0\end{array} \right)$.
As for $(\alpha, S_{E-v})$, we always have $\rho(\alpha, S_{E-v})\in [0,\frac{1}{2}]$.

 Let $\mu_{\alpha,v,\phi}$ be the spectral measure of $H_{ v,\alpha, \phi}$, and the \textit{integrated density of states}  $N_{\alpha,v}:\R\rightarrow [0,1]$  is defined as
\[N_{\alpha,v}(E)=\int_{\T^d}\mu_{\alpha,v,\phi}(-\infty,E]d\phi.\]
There is a close relation between  $\rho(\alpha, S_{E-v})$ and  $N_{\alpha,v}(E)$ \cite{AJ1}:
  \[N_{\alpha,v}(E)=1-2\rho(\alpha, S_{E-v}).\]

For any $g\in C^0(\T^d,\R)$, denote $R_{g}:=\left(
\begin{array}{ccc}
 \cos 2\pi g &   -\sin  2\pi g\cr
 \sin  2\pi g & \cos  2\pi g\end{array}
   \right)=e^{-2\pi g J}   ,$
  where $J=\left(\begin{array}{cc}0 & 1\\ -1 & 0\end{array}\right)$. Then we have the following
 facts about fibered rotation number (one can refer to \cite{Am09} for more details).

 \begin{Lemma}\label{Lem-rot-num}
 Let  $\alpha\in \T^d$ be rational independent and $A\in C^0(\T^d,SL(2,\R))$ be homotopic to identity.
 We have: \\
 (1) For rotation $R_{\varrho}$ with $\varrho\in \R$, there is
 \[|\rho(\alpha,A)-\varrho|\leq \|A-R_{\varrho}\|_0.\]
 (2) For general $C\in SL(2,\R)$, there exists  $c=c(C)\geq1$, such that
 \[|\rho(\alpha,A)-\rho(\alpha,C)|\leq c\|A-C\|_0^{1/2}.\]
\end{Lemma}

Usually the fibered rotation number varies after a  conjugation via some $B\in C^0(\T^d,PSL(2,\R))$, but the change depends only on the \textit{topological degree} of $B$. We say that the topological degree of $B$ is $\textbf{d}\in\Z^d$, as long as
$B$ is homotopic to
\[
\left(
\begin{array}{ccc}
 \cos \pi\langle \textbf{d}, \theta\rangle &   -\sin\pi\langle \textbf{d}, \theta\rangle\cr
 \sin \pi\langle \textbf{d}, \theta\rangle & \cos \pi\langle \textbf{d}, \theta\rangle\end{array}
   \right),\]
and we use $deg(B)$  to denote it.

 \begin{Lemma}[\cite{Kr}]\label{rota-number-rem}
 Let $A\in C^0(\T^d,SL(2,\R))$ and it is homotopic to identity, $B\in C^0(\T^d,PSL(2,\R))$. If $\widetilde{A}\in C^0(\T^d,SL(2,\R))$  and
 $(\alpha, A)$ is conjugate to $(\alpha, \widetilde{A})$
via $B$, then $\widetilde{A}$ is also homotopic to identity and
\[\rho(\alpha,\widetilde{A})=\rho(\alpha,A)+\frac{1}{2}\langle deg(B),\alpha \rangle \quad   (mod \quad \Z). \]
In particular, $\rho(\alpha,\widetilde{A})=\rho(\alpha,A)$  if $B$ is homotopic to identity.
 \end{Lemma}

\section{Proof of main results}\label{sec-2-freq-construction}

The proof  relies on the following well-known criterion of ac spectrum, which comes from the subordinacy theory.

\begin{Proposition}[{\cite{GP}}]\label{Lem-sub}
Let $\widetilde{\mathcal O}$ be the set of $E\in\R$ such that the cocycle $(\alpha, S_{E-v})$ is bounded.
 If $\widetilde{\mathcal O}$ has positive
Lebesgue measure, then  $H_{v,\alpha,\phi} $ has absolutely continuous
spectrum for any $\phi$.
\end{Proposition}

To obtain energies $E$, such that the corresponding cocycle $(\alpha, S_{E-v})$ is bounded, our proof is purely dynamical. Recall that rotations reducible cocycles are always bounded. Then the key dynamical results are the following:

\begin{Theorem}\label{thm-ana-tech}
Let $\varrho_0\in\R$,  $\mathfrak r>0$,   $F\in C^\omega_{\mathfrak r}(\T^d,sl(2,\R))$.     Suppose that  $\alpha\in WL$,   and the fibered rotation number
 $\rho(\alpha, R_{\varrho_0}e^{F})=:\varrho\in D_\alpha(\gamma',\tau')$, i.e. there exist $\gamma'>0,\tau'>0$ such that
   \begin{equation}\label{diophan-2}
\|2\varrho+\langle k, \alpha\ra\|_{\R/\Z}\geq \frac{\gamma'}{(1+|k|)^{\tau'}},\quad \forall k\in\Z^d.
\end{equation}
Then there exists $\varepsilon_*=\varepsilon_*(\alpha, \gamma', \tau',\mathfrak r)>0$, such that
if   $\|F\|_{\mathfrak r}<\varepsilon_*$, then
   $(\alpha,R_{\varrho_0}e^F)$ is $C^\infty$ rotations reducible.
\end{Theorem}

This result is for the multifrequency case, and if the frequency is two dimensional, then we first need to construct the super-Liouvillean frequency.
The construction of the frequencies is as follows: for any two real numbers $\tilde\alpha, \alpha'\in (0,1)$, suppose that $\{\tilde q_n\}_{n\in\N}, \{ q_n'\}_{n\in\N}$ are the convergences of $\tilde\alpha$ and $ \alpha'$ respectively.

\begin{Lemma}\label{freconstruct}
For any $\chi\geq 5$,  denote $\Omega(\chi)\subseteq \T^2$ the set of rational independent $\alpha=(\tilde\alpha,  \alpha')\in \T^2$  satisfying the following properties:
\begin{enumerate}
\item[$(a)$] \ $\tilde q_n> e^{ q_{n-1}'}$;
\item[$(b)$] \ $\tilde q_n^5< q_n'< 4\tilde q_n^\chi$;
\item[$(c)$] \ $\tilde q_n\wedge  q_{n-1}'=1$;
\item [$(d)$]\  $q_n'\wedge \tilde q_n=1$.
\end{enumerate}
Then $
  \Omega(\chi) \neq \varnothing$ and $ \Omega(\chi)\subseteq\{  \alpha \in\T^2\big|\ \beta(\alpha)=\infty\}.
$
\end{Lemma}

\begin{pf}
 Let us first prove that $\Omega(\chi)\neq \varnothing$. We construct the sequences $\{\tilde a_n\}_{n\geq 1}, \{ a_n'\}_{n\geq 1}$ of $\tilde\alpha, \alpha'$ respectively by induction. For $n=1$, we let $\tilde a_1=2$ and $a_1'=37$. Then $q_1'$ and $\tilde q_1$ satisfy (a)-(d) for $n=1$. For $n\geq 2$, suppose that we have already chosen $\tilde a_j,   a_j', j\leq n-1$, such that ($a$)-($d$) hold until the $n-1^{th}$ step. Since $\tilde q_{n-1}$ and $q_{n-1}'$ are coprime, there exist  $\tilde s_n, \tilde t_n\in\Z$ with $|\tilde s_n|<\tilde q_{n-1}\tilde q_{n-2}, |\tilde t_n|<q_{n-1}'\tilde q_{n-2}$, such that
\[\tilde t_n\tilde q_{n-1}+\tilde q_{n-2}=\tilde s_n  q_{n-1}'.\]
Moreover, there exists $\tilde k_n\in\Z_+$, such that $\tilde k_n\wedge  q_{n-1}'=1$, and
\[\tilde k_n\tilde q_{n-1}\geq 2 e^{ q_{n-1}'}. \]
Let $\tilde a_n=\tilde t_n+\tilde k_n\in\Z_+$. Then we have
\[\tilde q_n=\tilde a_n\tilde q_{n-1}+\tilde q_{n-2}=\tilde k_n\tilde q_{n-1}+\tilde t_n\tilde q_{n-1}+\tilde q_{n-2}.\]
Thus, we obtain that $\tilde q_n\wedge  q_{n-1}'=1$, and
\[\tilde q_n\geq 2e^{ q_{n-1}'}-  q_{n-1}'\tilde q_{n-1}\tilde q_{n-2}>  e^{ q_{n-1}'},\]
since $q_{n-1}'\geq q_1'=37$.
As for the construction of $ a_n'$, we follow the same idea: let
$  a_n'= t_n'+  k_n'$, where $ t_n'  q_{n-1}'+q_{n-2}'$ divides $\tilde q_n$,  $  k_n'\wedge \tilde q_n=1$, and
$2\tilde q_{n}^5\leq   k_n'  q_{n-1}'\leq  3\tilde q_n^\chi.$
Then we get $(b)$ and $(d)$ at  the $n^{th}$ step. We thus complete the proof of $\Omega(\chi)\neq \varnothing$.

Moreover, for any $\alpha\in\Omega(\chi)$, by properties (a) and (b), we get for any $n\geq 1$ that $\tilde q_n>e^{\tilde q_{n-1}^5}$. Therefore, we have $ \Omega(\chi)\subseteq\{  \alpha \in\T^2\big|\ \beta(\alpha)=\infty\}$.

\end{pf}

\begin{Remark}\label{compareaj}
In Avila-Jitomirskaya's result \cite{AJ2}, they choose the alternating super-Liouvillean frequency, i.e.
\begin{enumerate}
\item[$(a')$] \ $  \tilde q_n > e^{q_{n-1}'}$;
\item[$(b')$] \ $q_n'> e^{\tilde q_{n}}$.
\end{enumerate}
Therefore, compared to our frequency, the main difference is $(b)$, while we can also choose $\tilde q_n$  arbitrarily
 larger than $q_{n-1}'$, we must assume $q_n'$  grows polynomially with respect to $\tilde q_{n}$.
\end{Remark}

Once we choose the frequency, then we have the following rotations reducibility result.

\begin{Theorem}\label{thm-2-fre-tech}
 Let $\mathfrak{r},\gamma', \tau'>0, \chi\geq 5$,    $F\in C^\omega_{\mathfrak r}(\T^d,sl(2,\R))$. Suppose that   $\alpha\in\Omega(\chi)$  and
 $\rho(\alpha, R_{\varrho_0}e^{F})\in D_\alpha(\gamma',\tau')$. Then there exists $\varepsilon_*=\varepsilon_*(\alpha, \gamma', \tau',\mathfrak r)>0$, such that
if   $\|F\|_{\mathfrak r}<\varepsilon_*$, then
   $(\alpha,R_{\varrho_0}e^F)$ is $C^\infty$ rotations reducible.
\end{Theorem}

The proof of Theorem \ref{thm-ana-tech} and  Theorem \ref{thm-2-fre-tech}  will be left to section \ref{sec4} and section \ref{sec:6}. Now let's finish the proof  of Theorem \ref{thm-ana} and Theorem \ref{thm-2-fre} based on  the above results.  Here we only give the proof of
Theorem \ref{thm-2-fre}, and the proof of Theorem \ref{thm-ana} is the same. In order to apply Proposition \ref{Lem-sub} to prove Theorem \ref{thm-2-fre}, we have to obtain positive Lebesgue measure of the energies $E$ such that the corresponding cocycles are bounded, while Theorem \ref{thm-2-fre-tech} only gives us positive Lebesgue measure of the fibered rotation numbers (not the energies) such that the corresponding cocycle is bounded. Hence, we need the following elementary result, and this kind of idea first appeared in \cite{ABD09}.

\begin{Lemma}[{\cite[Lemma~3.9]{BK}}]\label{Lem-Lip-positive}
Let $I, J$ be two intervals of $\R$ and $f: I\rightarrow J$ be a continuous map. Assume that there exists a set $\Upsilon\subseteq J$, of positive Lebesgue measure, such that $\forall\, x\in f^{-1}(\Upsilon)$ there is a constant $c_x>0$ such that $|f(x)-f(x')|\leq c_x |x-x'|$ for all $x'$ close to $x$. Then $ f^{-1}(\Upsilon)$ has positive Lebesgue measure.
\end{Lemma}

Once we have this, we can now finish the whole proof. \\

\textbf{Proof  of
Theorem \ref{thm-2-fre}.}  By Lemma \ref{freconstruct}, we can select $\alpha \in
  \Omega(\chi)  \subseteq \{  \alpha \in\T^2\big|\ \beta(\alpha)=\infty\}
$.  Fix some absolutely small $0<\sigma<<1$. Note
$S_{E- \lambda v}=\left(
\begin{array}{ccc}
 E &   -1\cr
   1 & 0\end{array}
   \right)\exp\{ \lambda v(\phi)\left(
\begin{array}{ccc}
 0 &   0\cr
   1 & 0\end{array}
   \right)\}$.
 One can find $E_1<E_2$, such that
\begin{equation}\label{rot-0-using-in-proof-of-thm-2-fre}
[2\sigma,\,\frac{1}{2}-2\sigma]\subseteq\rho_0([E_1,\,E_2])\subseteq[\sigma,\,\frac{1}{2}-\sigma],
\end{equation}
where we denote $\rho_0(E)=\rho(\alpha,\left(
\begin{array}{ccc}
 E &   -1\cr
   1 & 0\end{array}
   \right))$.
Then, $\forall \, E\in [E_1,E_2]$, $\exists \, P_E\in SL(2,\R)$,
such that $P_E\left(
\begin{array}{ccc}
 E &   -1\cr
   1 & 0\end{array}
   \right)P_E^{-1}=R_{\rho_0(E)}$.
Moreover, $P_E$ are uniformly bounded on $[E_1,E_2]$.

Denote $\rho_{\lambda v}(E)=\rho(\alpha,S_{E-\lambda v})$ for short.
By $(\ref{rot-0-using-in-proof-of-thm-2-fre})$ and Lemma \ref{Lem-rot-num}, there exists $\epsilon'>0$, when $\|\lambda v\|_{\mathfrak r}<\epsilon'$,  we have $|\rho_{\lambda v}(E)-\rho_0(E)|<\sigma$,
which implies that $\rho_{\lambda v}([E_1,\,E_2])\supseteq [3\sigma,\,\frac{1}{2}-3\sigma]$ and then
$\rho_{\lambda v}([E_1,\,E_2])\cap D_{\alpha}(\gamma',\tau')$ is of positive measure. Then $P_E$ conjugates $(\alpha,\, S_{E-\lambda v})$ to
$(\alpha,\, R_{\rho_0(E)}e^{F_E})$, where $F_E(\cdot):= \lambda v(\cdot)P_E\left(
\begin{array}{ccc}
 0 &   0\cr
   1 & 0\end{array}
   \right)P_E^{-1}$. Let $\varepsilon_*>0$ be the one given in Theorem \ref{thm-2-fre-tech}.  There exists $\epsilon''\in(0,\epsilon']$, such that when $\|\lambda v\|_{\mathfrak r}<\epsilon''$ and $\rho_{\lambda v}(E)\in  D_\alpha(\gamma',\tau')$, we have $\|F_E\|_{\mathfrak r}<\varepsilon_*$ and  $\rho(\alpha, R_{\rho_0(E)}e^{F_E})\in  D_\alpha(\gamma',\tau')$. Therefore, by Theorem \ref{thm-2-fre-tech},
$(\alpha,\, R_{\rho_0(E)}e^{F_E})$  is rotations reducible and so is $(\alpha,\, S_{E-\lambda v})$.
Let $\Pi(\gamma',\tau'):=\{E\in [E_1,E_2]\,\big|\, \rho_{\lambda v}(E)\in D_{\alpha}(\gamma',\tau')\}$. When $\|\lambda v\|_{\mathfrak r}<\epsilon''$ and $E\in  \Pi(\gamma',\tau')$,
then $(\alpha,\, S_{E-\lambda v})$  is  rotations reducible.  In order to apply Lemma \ref{Lem-Lip-positive}, we need the following lemmas:

 \begin{Lemma}\label{Lem-rot-num-1}
 Let  $\alpha\in \T^d$ be rational independent and $A\in C^0(\T^d,SL(2,\R))$ be homotopic to identity.
For any $SO(2,\R)$-valued $R_g$ with $g\in C^0(\T^d,\R)$, there exists an absolute constant $c_*\geq 1$, such that if $\|A-R_g\|_{0}<\frac{1}{4}$, then
 \[|\rho(\alpha,A)-\rho(\alpha,R_g)|\leq c_*\|A-R_g\|_0.\]
\end{Lemma}
\begin{pf}
 Denote $\tilde A:=R_{-g}A$, and then $A=R_g\tilde A$. Then $(\alpha, A)=(\alpha, R_g)\circ (0,\tilde A)$, and thus
\[F_{(\alpha,A)}(\phi,x)=F_{(\alpha, R_g)}\circ F_{(0,\tilde A)}(\phi,x)=(\phi+\alpha, x+g(\phi)+d_{\tilde A}(\phi,x)),\]
where $F_{(\alpha, A)}$ is a lift of $f_{(\alpha, A)}$ as in section \ref{sec-cocycles}.
Since $F_{(\alpha, R_g)}(\phi,x)=(\phi+\alpha, x+g(\phi))$, then by the definition of fibered rotation number, we get that
\[|\rho(\alpha, A)-\rho(\alpha, R_g)|\leq \|d_{\tilde A}\|_{0}.\]
By a direct computation, there exists $c_*\geq 1$ such that  for $\|A-R_g\|_{0}<\frac{1}{4}$ we have
\[
\|d_{\tilde A}\|_{0}\leq c_*\|\tilde A-I\|_{0}\leq c_*\|R_{-g}\|_{0}\|A-R_g\|_{0}\leq c_*\|A-R_g\|_{0}.\]
Consequently, we have
 \[|\rho(\alpha, A)-\rho(\alpha, R_g)|\leq c_*\|A-R_g\|_{0}.\]
 \end{pf}

\begin{Lemma}\label{Lem-rot-to-Lip}
If $(\alpha, S_{E-v})$ is rotations reducible, there exists a constant $c_E>0$, such that
 $|\rho_v(E)-\rho_v(E')|\leq c_E |E-E'| $   for all $E'$ close to $E$.
\end{Lemma}

\begin{pf}
We assume that $(\alpha, S_{E-v})$ is conjugate via $B$ to $(\alpha, R_{g})$. Since $(\alpha, S_{E-v})$ is homotopic to the identity, the same is true for $(\alpha, R_g)$ by Lemma \ref{rota-number-rem}. Thus, we have $g\in C^0(\T^2,\R)$.  Moreover, by Lemma \ref{rota-number-rem}, we obtain that $\rho_v(E)=\rho(\alpha, R_g)-\frac{1}{2}\la deg(B), \alpha\ra (\textrm{mod}\ \Z)$ and $\rho_v(E')=\rho(\alpha, B(\cdot+\alpha)S_{E'-v}(\cdot)B^{-1}(\cdot))-\frac{1}{2}\la deg(B), \alpha\ra (\textrm{mod}\ \Z)$.
Let $F(\cdot)=B(\cdot+\alpha)\left(
\begin{array}{ccc}
 1 &   0\cr
   0 & 0\end{array}
   \right)B^{-1}(\cdot)$. Then obviously, we have $\|F\|_{0}\leq \|B\|_{0}^2$. For $E'=E+\Delta E$,  there is
\begin{eqnarray*}
 B(\cdot+\alpha)S_{E'-v}(\cdot)B^{-1}(\cdot)=R_g(\cdot)+\Delta E \cdot F(\cdot).
\end{eqnarray*}
Then for $E'$ close to E, we have $\|\Delta EF\|_{0}<\frac{1}{4}$, and by Lemma \ref{Lem-rot-num-1}, we get
\[|\rho(\alpha, R_g)-\rho(\alpha, R_g+\Delta E \cdot F)|\leq c_*|\Delta E|\|F\|_{0}<c_E|\Delta E|,\]
which implies the desired conclusion.
\end{pf}

Consequently,  by Lemma \ref{Lem-rot-to-Lip} and Lemma \ref{Lem-Lip-positive},  $\Pi(\gamma',\tau')$ is of positive measure. By Proposition  \ref{Lem-sub}, we thus finish the whole proof. \qed

\section{Proof of Theorem \ref{thm-ana-tech}} \label{sec4}
In this section, we use KAM iteration to prove   Theorem \ref{thm-ana-tech}.  The philosophy of KAM iteration is to find a sequence of coordinate transformations through  which the perturbations gradually become smaller and smaller. In the end, the transformations converge and  the limit conjugates the original system to some simpler one.

However, our proof relies on a modified KAM scheme.  In conjugating  a cocycle $(\alpha, R_\varrho e^{F})$ to some simpler one, a crucial ingredient is to solve the linearized homological  equations in the form
\begin{eqnarray*}
h_1(\cdot+\alpha)-h_1+f_1=f_1^{(u)},\quad e^{4\pi i\varrho}h_2(\cdot+\alpha)-h_2+f_2=0,
\end{eqnarray*}
 where $h_1,h_2$ are unknown, $f_1,f_2$ comes from $F$, and $f_1^{(u)}$ is the unsolvable part. When $\alpha$ is Diophantine, $f_1^{(u)}$ is a constant and $(\alpha, R_\varrho e^{F})$  is conjugated to $(\alpha, R_{\varrho_+} e^{F_+})$ ($\varrho_+\in\R$) with  $F_+$ being much smaller. We then solve new linearized homological  equations in similar form \cite{E92}. But, for $\alpha\in WL$, $f_1^{(u)}$ will depend on $\varphi:=\phi_1\in\T^1$, where $\phi=(\phi_1,\cdots,\phi_d)\in\T^d$. As a result, we arrive at $(\alpha, R_{\varrho+\frac{g(\varphi)}{2\pi}}e^{F_+(\phi)})$ with  $F_+$ being much smaller while $g$ is as large as $F$. Thus, in our KAM iteration, we  should consider linearized homological  equation in the form
 \begin{eqnarray*}
e^{2i(2\pi \varrho+g(\varphi))}h(\cdot+\alpha)-h+f=0,
\end{eqnarray*}
with unknown $h$. It is the essential difference from the classical one \cite{E92} and we deduce it using techniques in \cite{KWYZ}.

In the sequel,  for $\phi=(\phi_1,\cdots,\phi_d)\in\T^d$, we denote $\varphi:=\phi_1\in\T^1, \theta:=(\phi_2, \cdots, \phi_d)\in\T^{d-1}$ for simplicity.

\subsection{Linearized Homological Equation}
To prove Theorem \ref{thm-ana-tech}, one needs  to consider  the  linearized homological equation of the form
\begin{equation}\label{homo-equ-f-pre}
e^{2i(2\pi\varrho+g(\varphi))}h(\phi+\alpha)-h(\phi)=-f(\phi)
\end{equation}
with $f \in C_{r,s}^\omega(\T^{1}\times \T^{d-1}, \C)$ and $g\in
C_{r}^\omega(\T^1,\R)$ for some $r,s>0$. Such an equation can be furthermore written in  the form
\begin{equation}\label{homo-equ-f}
(e^{4\pi i\varrho}+\widetilde{g}(\varphi))h(\phi+\alpha)-h(\phi)=-f(\phi)
\end{equation}
where $\widetilde{g}\in C_{r}^\omega(\T^1,\R)$ defined as
\begin{eqnarray}\label{equ-tilde-g}
\widetilde{g}(\varphi)=e^{4\pi i\varrho}(e^{2i g(\varphi)}-1).
\end{eqnarray}
Solving such kind of  linearized homological equations will play an important role in  the KAM-type proof of
Theorem \ref{thm-ana-tech}.  We solve (\ref{homo-equ-f}) by an approximate equation, which is the content of the following proposition.

\begin{Proposition}\label{homo-lemma}
Let $\gamma'>0, \tau'>0, r>\sigma>0, s>\delta>0, \sigma\leq
\delta\leq \frac{1}{4}$, $0<\tilde \eta\leq \eta \ll 1$, and $g\in
C_{r}^\omega(\T^1,\R)$, $f\in C_{r,s}^\omega(\T^1\times \T^{d-1},\C) $.
If
$\|g\|_{r}\leq \eta$, $\|f\|_{r,s}\leq
\tilde \eta$,
\begin{equation}\label{dia-2}
K=\left[\frac{1}{2\pi\sigma}\ln \frac{1}{\tilde \eta}\right]+1<
\left(\frac{\gamma' \sigma}{32\pi\eta}\right)^{\frac{1}{ \tau' }},
\end{equation}
and $\varrho\in D_\alpha(\gamma',\tau')$,
then  (\ref{homo-equ-f-pre}) has an approximate
solution $h\in C_{r-\sigma, s-\delta}^\omega(\T^1\times \T^{d-1}, \C)$
with
$$\|h\|_{r-\sigma, s-\delta}\leq \frac{C_0(d)K^{\tau'} }{\gamma'\sigma^{d}}\cdot \tilde\eta,$$
and the error term $\tilde
P=\mathcal{R}_K( \tilde g(\varphi)h(\phi+\alpha
)+f(\phi))$
with $\tilde g$ satisfying (\ref{equ-tilde-g}) and
$$\|\tilde P\|_{r-2\sigma, s-2\delta}\leq C_1(d)K^{d}\left(1+\frac{K^{\tau'} }{\gamma'\sigma^{ d}}\eta\right)\cdot \tilde \eta^2,$$
where $C_0,C_1$ are constants only depending on $d$.
\end{Proposition}

The proof of Proposition \ref{homo-lemma} is somehow technical and we put it in the Appendix.

\subsection{KAM Step}\label{sec-kam-conti}
In this section, we will give the inductive lemma and details about
one step of the iteration.

Before stating the inductive lemma, we first give some notations. Let $\alpha :=(\tilde\alpha, \alpha')$.
For any $r,s,\eta,\tilde \eta,\gamma',\tau'>0$, we define
\begin{eqnarray*}
 \lefteqn{\mathcal{F}_{r,s}(\eta, \tilde{\eta},\gamma',\tau'):=}\\&&\left\{
\begin{array}{cc}
(\alpha, \,R_{\varrho+\frac{g(\varphi)}{2\pi} }e^{F(\varphi,\theta)})\in
C^\omega_{r,s}(\T^1 \times \T^{d-1},SL(2,\R)):\\
\rho(\alpha, R_{\varrho+\frac{g}{2\pi}} e^F)\in D_\alpha(\gamma',\tau') \ \
\textrm{and}\  \|g\|_{r} \leq \eta, \|F\|_{r,s}\leq
\tilde{\eta}
\end{array}
\right\}.
\end{eqnarray*}

Without loss of generality,
we assume that there is some  $\gamma,\tau>0$ and $\tilde U\in [0,+\infty)$, such that $\alpha=(\tilde\alpha,\alpha') \in WL(\gamma,\tau,\tilde U)$. Suppose that  $(Q_j)$ is the
selected sequence of $\tilde\alpha$ by Lemma \ref{cd-coroll} with
$\mathcal A=4$ and $U=\tilde U(\tilde \alpha)+8$.
For $\gamma, \gamma'>0, \tau, \tau'>0$, let
 \[\tau_*=\max\{\tau,\tau'\},\ \  \gamma_*=\min\{\gamma, \gamma'\}.\]
 Let $r_0>0, s_0>0$, and $Q_*\in\N$ be the smallest integer such that for any $Q> Q_*$ we have
\[\ln Q<\frac{Q^{1/8}r_0}{40c \tau_* U},\]
where c is a global constant with $c>100$.
Now suppose $\varepsilon_0$ is small enough such that
\begin{equation}\label{pert-condi}
\varepsilon_0<\min\{\frac{(r_0s_0\gamma_*C_*^{-1})^{60(\tau_*+d)}}{(\tau_*+d)!Q_1^{12c\tau_*U}},
e^{-2c\tau_*U},    e^{-40(\ln Q_*)^2c\tau_* U} \}, \ \
\ln\frac{1}{\varepsilon_0}<(\frac{1}{\varepsilon_0})^{\frac{1}{12\tau_*}}.
\end{equation}
For any
given $r_0,s_0,\varepsilon_0$ satisfying (\ref{pert-condi}),  let $\Delta_0=s_0/4$, and we
inductively define some sequences depending on
$r_0,s_0,  \Delta_0,\varepsilon_0$ for $j\geq 1$:
\begin{equation}\label{para-induc}
\left.\begin{array}{ll}  r_j=\frac{r_0}{2Q_j^3} \ \ &
 \bar r_j=\frac{r_0}{Q_j^3},  \\   \Delta_j=\Delta_0/2^j, \ \ & s_j=s_{j-1}-\Delta_{j},\\
\varepsilon_j=\frac{\varepsilon_{j-1}}{Q_{j+1}^{2^{j+1}c\tau_*U}}, \
\ &
\tilde\varepsilon_j=\sum_{m=0}^{j-1}\varepsilon_m.
\end{array}
\right.
\end{equation}

\begin{Proposition}[Iterative Lemma]\label{iter-lemma}
For any $\varepsilon_0, r_0, s_0, \gamma, \gamma', \tau, \tau'>0$ satisfying (\ref{pert-condi}),   and $\alpha=(\tilde\alpha,\alpha') \in WL(\gamma, \tau,\tilde U)$ with $\tilde U=\tilde U(\tilde\alpha)<\infty$, we define
$\varepsilon_n, \tilde\varepsilon_n, r_n, s_n$ as in
(\ref{para-induc}). Then the following holds
for $n\geq 2$: If the cocycle
\begin{equation}\label{sys-iter-former}
(\alpha,\, R_{ \varrho+\frac{g_{n-1}(\varphi)}{2\pi}}e^{F_{n-1}(\varphi,\theta)})
 \in \mathcal{F}_{r_{n-1},s_{n-1}}(4\tilde\varepsilon_{n-1},\varepsilon_{n-1},\gamma', \tau'),
\end{equation}
then it is conjugate to the cocycle
\begin{equation}\label{sys-iter-after}
(\alpha,\,R_{ \varrho+\frac{g_{n}(\varphi)}{2\pi}}e^{F_{n}(\varphi,\theta) })
\in \mathcal{F}_{r_{n},s_{n}}(4\tilde\varepsilon_{n},\varepsilon_{n},\gamma', \tau'),
\end{equation}
by an analytic conjugation $\Phi_n\in C_{ r_{n},s_{n}}^\omega(\T^1\times
\T^{d-1}, SL(2,\R))$ satisfying $$\|\Phi_n-I\|_{r_{n},s_{n}}\leq
2\varepsilon_{n-1}^{3/4}.$$
\end{Proposition}\medskip


We will divide the proof of Proposition \ref{iter-lemma} into
different lemmas.

\begin{Lemma}\label{reson-lemma}
 For $n\geq 2$, any cocycle
\begin{equation}\label{sys-origin}
(\alpha,\,R_{ \varrho+\frac{g(\varphi)}{2\pi}}e^{F(\varphi,\theta)})
 \in \mathcal F_{r_{n-1}, s_{n-1}}(4\tilde\varepsilon_{n-1},\varepsilon_{n-1},\gamma', \tau')
\end{equation}
can be conjugate to another cocycle
\begin{equation}\label{sys-after-reso}
(\alpha,\,R_{ \rho_f+\frac{\tilde g(\varphi) }{2\pi}}e^{\tilde F(\varphi,\theta)})
 \in \mathcal F_{\bar r_n, s_{n-1}}(\varepsilon_{n-1}^{3/4},2\varepsilon_{n-1}, \gamma', \tau'),
\end{equation}
where $\rho_f=\rho(\alpha, R_{ \varrho+\frac{g}{2\pi}}e^F)$
via some conjugation  $R_{\frac{v_n(\varphi)}{2\pi}}=e^{-v_n(\varphi)J}$, with $v_n$ satisfying
\[\|v_n\|_{\bar r_n}<Q_n^{\frac{7}{4}}\varepsilon_0^{1/2}.\]
Moreover, the fibered rotation number remains unchanged.\end{Lemma}
\begin{pf}
Let $v_n$ be the solution of
\[  v_n(\varphi+\tilde \alpha)-v_n(\varphi)=
-\mathcal T_{Q_n} g(\varphi)+\hat g(0).\]
 Then we have
$$\|v_n\|_{\frac{r_{n-1}}{2}}\leq 2Q_n\sum_{0<|k|<Q_n}\|g\|_{r_{n-1}}e^{-2\pi|k|r_{n-1}/2}\leq \frac{64Q_n\varepsilon_0 }{r_{n-1}}
\leq Q_n^{\frac{7}{4}}\varepsilon_0^{1/2},$$ and the given cocycle is conjugate to
$(\alpha,\,R_{\tilde \varrho+\frac{\mathcal R_{Q_n}g(\varphi)}{2\pi} }e^{\tilde F(\varphi,\theta)})$,
where
\[\tilde\varrho=\varrho+\frac{\hat g(0)}{2\pi},\quad
\tilde F(\varphi,\theta) = e^{-v_n(\varphi)J}F(\varphi,\theta)e^{v_n(\varphi)J}. \]
Here we
need a small trick used in \cite{YZ2,KWYZ}, saying $|\Im v_n(\varphi)|$ can be
well controlled by sacrificing the analytic radius. We give a short review of the proof for completeness.

\begin{Lemma}\label{claim-radius}
Let $r>0, \tilde\omega\in\R^m (m\in\N)$, and $v(\varphi)$ is the solution of $v(\varphi+\tilde\omega)-v(\varphi)=\mathcal T_K g(\varphi)-\hat g(0)$, where $g\in C_{r}^\omega(\T^m,\R)$ and $K\in \N$.  If $\|\la k,\tilde\omega\ra\|_{\R/\Z}\geq \iota^{-1}>0$ for all $0<|k|<K\ (k\in\Z^m)$, then for any $\bar r<r$, we have
$$\|\Im v \|_{\bar r}\leq\frac{C\bar r \iota\|g\|_{r}}{(r-\bar r)^{m+1}},$$
where $C$ is a global constant depending on $m$.
\end{Lemma}
\begin{pf}
Let $\varphi=\varphi_1+i\varphi_2$, with $\varphi_1\in \T^m, \varphi_2\in\R^m$. Denote
\[v_1(\varphi)=v(\varphi_1),\ \ \textrm{and}\ \ v_2(\varphi)=v(\varphi)-v_1(\varphi).\]
Owing to the fact that $g(\varphi)$ is real-analytic in $\varphi$, we get $v_1(\varphi)\in\R$. Therefore,
\begin{eqnarray*}
\lefteqn{\|\Im v \|_{\bar r}=\|\Im v_2 \|_{\bar r}}\\
&\leq& \sup_{ |\varphi_2|<\bar r} \sum_{0<|k|<K}\iota |\hat g(k)|\cdot |e^{-2\pi\la k,\varphi_2\ra}-1|\\
&\leq& \iota\sum_{0<|k|<K}\|g\|_{r} e^{-2\pi|k|(r-\bar r)}\cdot 2\pi|k|\bar r\leq \frac{C\bar r \iota\|g\|_{r}}{(r-\bar r)^{m+1}}.
\end{eqnarray*}
\end{pf}

By Lemma \ref{claim-radius}, we get (as $\varepsilon_0$ is sufficiently small)
\begin{equation}\label{esti-v-n}
\|\Im v_n\|_{\bar r_n} \leq \frac{C
r_0^{-1} \tilde \varepsilon_{n-1}}{Q_n^{1/2}} \ll 1.
\end{equation}
Thus, we have
\[\|\tilde F \|_{\bar r_n, s_{n-1}}<2\varepsilon_{n-1} .\]
 Moreover, standard computation shows that  for $\varepsilon_0$ small enough we have
$$\|\mathcal{R}_{Q_n}g\|_{ r_{n-1}/2}\leq CQ_ne^{-2\pi Q_nr_{n-1}/2}\|g\|_{r_{n-1}}
\leq
\varepsilon_{n-1}^{4/5},$$
by (\ref{pert-condi}) and the selection of $Q_*$.  Since $e^{v_n(\varphi)J}$ is homotopic to the
identity,  by Lemma \ref{Lem-rot-num} and \ref{rota-number-rem}, the fibered
rotation number of the new system remains $\rho_f$, and $|\rho_f-\tilde \varrho|\leq 2(\|\mathcal R_{Q_n}g\|_0+\|\tilde F\|_0)\leq 4\varepsilon_{n-1}^{4/5}$. Let $\tilde g(\varphi)=\mathcal R_{Q_n}g(\varphi)+2\pi(\tilde\varrho-\rho_f)$. Then
$$\|\tilde g\|_{\bar r_n}\leq \|\mathcal R_{Q_n}g\|_{r_{n-1}/2}+2\pi|\rho_f-\tilde\varrho|<\varepsilon_{n-1}^{3/4}.$$
We finish the proof of Lemma \ref{reson-lemma}.
\end{pf}\medskip

\begin{Remark}\label{large-remark}
In the proof of Lemma \ref{reson-lemma}, although the norm of the
transformation $e^{v_nJ}$ is not large, it may  not be close to the
identity.
\end{Remark}

Once we get the cocycle (\ref{sys-after-reso}), we will use Proposition \ref{homo-lemma} to make
the perturbation small enough.

\begin{Lemma}\label{redu-lemma}
Under the assumptions of Lemma \ref{reson-lemma}, for $n\geq 1$, 
there exists a conjugation map $\Psi_n\in
C_{r_n, s_{n}}^{\omega}(\T^1\times \T^{d-1}, SL(2,\R))$ with
\begin{equation}
\|\Psi_n- I\|_{ r_n, s_n}\leq \varepsilon_{n-1}^{3/4},
\end{equation}
such that $\Psi_n$ conjugates (\ref{sys-after-reso}) to
\begin{equation}\label{sys-after-iter}
(\alpha,\,R_{ \rho_f+\frac{\tilde g_+(\varphi)}{2\pi} }e^{\tilde F_+(\varphi,\theta)}) \in \mathcal F_{r_n, s_n}( 2\varepsilon_{n-1}^{3/4},
\varepsilon_n/2,\gamma', \tau'),
\end{equation}
satisfying
\[\|\tilde g_+-\tilde g\|_{r_n}\leq 4\varepsilon_{n-1}.\]

\end{Lemma}
\begin{pf}
For simplicity, we denote temporarily $\tilde r_0:=\left\{\begin{array}{ll} r_0, & n=1 \\ \bar r_n=\frac{r_0}{Q_n^3}, & n\geq 2
\end{array},\right.$
$\tilde s_0:=s_{n-1},$ $\tilde\Delta:=\Delta_n, \eta:=2\varepsilon_{n-1}^{3/4}, \tilde
\eta_0:=2\varepsilon_{n-1}, \tilde g_0(\varphi):=\tilde g(\varphi), \tilde
F_0(\varphi,\theta):=\tilde F(\varphi,\theta)$, and define sequences
inductively for $\nu\geq 1$: \footnote{Note that $\frac{1}{Q_n^3}$ goes to 0 much
faster than $\Delta_n$, we can assume $\tilde r_0<2\tilde\Delta$ in
general, and thus $\sigma_\nu< \delta_\nu$.}
\begin{displaymath}\left.
\begin{array}{ll}
\sigma_{1}=\frac{\tilde{r_0}}{8},&
\sigma_{\nu}=\frac{1}{2^{\nu-1}}\sigma_{1},\\
\delta_{1}=\frac{\tilde\Delta}{4},&
\delta_{\nu}=\frac{1}{2^{\nu-1}}\delta_{1},\\
\widetilde{r}_{\nu}=\widetilde{r}_{\nu-1}-2\sigma_{\nu},&
\widetilde{s}_{\nu}=\widetilde{s}_{\nu-1}-2\delta_{\nu},\\
\tilde\eta_{\nu}=\tilde\eta_0^{(\frac{3}{2})^{\nu}}, &
K_{\nu}=[\frac{1}{2\pi\sigma_{\nu}}\ln\frac{1}{\tilde\eta_{\nu-1}}]+1.
\end{array}
\right.
\end{displaymath}
Let $N=[2^nc_1\tau_*U\ln Q_n]+1$, where
$c_1=\frac{c}{24\tau_*\ln 3}$. For convenience, we denote \[\tilde A_{j}= R_{\rho_f+\frac{\tilde
g_{j}}{2\pi}}.\] Assume that for $j=1,2,\cdots,
\nu-1< N\ (\nu\geq 1)$, there exist $Y_j, \tilde F_j\in C_{\tilde
r_j,\tilde s_j}^\omega(\T^1\times \T^{d-1}, sl(2,\R))$,  with $\|Y_j\|_{\tilde r_j,\tilde
s_j}\leq \tilde\eta_{j-1}^{4/5}$, 
$\|\tilde F_j\|_{\tilde r_j,\tilde s_j}\leq \tilde
\eta_j$, so that $e^{Y_j}$ conjugates the cocycle
$(\alpha,\,\tilde A_{j-1}e^{\tilde F_{j-1}(\varphi,\theta)})$ to the cocycle
$(\alpha,\,\tilde A_j e^{\tilde F_{j}(\varphi,\theta)})$,
where $\tilde F_j=\left(
\begin{array}{ccc}
\tilde  F_j^{11} & \tilde F_j^{12}\cr \tilde F_j^{21} & -\tilde
F_j^{11}\end{array} \right)$ and $\tilde g_j=\tilde
g_{j-1}+\frac{ [\tilde F_{j-1}^{12}-\tilde
F_{j-1}^{21}]_\theta}{2} $.

The linearized homological equation in  $\nu-1$-th step is
\begin{equation}\label{homo-matri}
\tilde A_{\nu-1}^{-1}Y_{\nu}(\cdot +\alpha)\tilde A_{\nu-1}- Y_{\nu}+\tilde F_{\nu-1}\equiv [\tilde F_{\nu-1}^-]_\theta J.
\end{equation}
Let $M=\frac{1}{1+i}\left(
\begin{array}{ccc}
  1 & -i\cr   1 & i\end{array} \right)\in U(2)$. Then 
$\tilde F_{\nu-1}$ can be uniquely
rewritten as
\[
\tilde F_{\nu-1}=\left(
\begin{array}{ccc}
  \tilde F_{\nu-1}^{11} &  \tilde F_{\nu-1}^{+}+\tilde F_{\nu-1}^-\cr  \tilde F_{\nu-1}^{+}-\tilde F_{\nu-1}^- & -
\tilde F_{\nu-1}^{11}\end{array} \right)=M^{-1}\left(
\begin{array}{ccc}
   i\tilde F_{\nu-1}^{-} &   w_{\nu-1}^1 \cr      w_{\nu-1}^2 & -i\tilde F_{\nu-1}^{-}\end{array} \right)M,\]
where $\tilde F_{\nu-1}^{\pm}=\frac{1}{2}(\tilde F_{\nu-1}^{12}\pm \tilde F_{\nu-1}^{21}), w_{\nu-1}^1=\tilde F_{\nu-1}^{11}-i\tilde F_{\nu-1}^+,
 w_{\nu-1}^2=\tilde
F_{\nu-1}^{11}+i\tilde F_{\nu-1}^+$. Similarly, we can rewrite
\[Y_{\nu}=M^{-1}\left(\begin{array}{ccc}iY_{\nu}^- & h_{\nu}^1\cr h_{\nu}^2 & -iY_{\nu}^- \end{array}\right)M,\]
 where
$Y_{\nu}^{\pm}=\frac{1}{2}(Y_{\nu}^{12}\pm Y_{\nu}^{21}), h_{\nu}^1= Y_{\nu}^{11}-i Y_{\nu}^+,
 h_{\nu}^2=
Y_{\nu}^{11}+i Y_{\nu}^+$. 
Then the homological equation (\ref{homo-matri}) is equivalent to
\begin{eqnarray}
\label{homo-matri-1}   Y_{\nu}^-(\cdot+\alpha)-Y_{\nu}^-=-\tilde F_{\nu-1}^-+[\tilde F_{\nu-1}^-]_\theta,\\
\label{homo-matri-2}e^{ 2i(2\pi\rho_f+\tilde g_{\nu-1})}h_{\nu}^1(\cdot+\alpha)-h_{\nu}^1=-w_{\nu-1}^1,\\
\label{homo-matri-3} e^{ -2i(2\pi\rho_f+\tilde g_{\nu-1})}h_{\nu}^2(\cdot+\alpha)-h_{\nu}^2=-w_{\nu-1}^2.
\end{eqnarray}
Since $(\tilde\alpha, \alpha')\in WL(\gamma,\tau,\tilde U)$,
it is classical to solve (\ref{homo-matri-1}) with
 \[\|Y_{\nu}^-\|_{\tilde r_{\nu},\tilde s_{\nu}}\leq
\frac{4(\tau+d)!\tilde \eta_{\nu-1}}{\gamma
\sigma_{\nu}^{\tau+d}}<\tilde \eta_{\nu-1}^{4/5}.\]
We will apply Proposition \ref{homo-lemma} to solve (\ref{homo-matri-2}).
First, we can check that
\begin{eqnarray*}
&&\|\tilde g_{\nu-1}\|_{ \tilde
r_{\nu-1}}\leq \|\tilde g_0\|_{\tilde r_0}+\sum_{j=0}^{\nu-2}\tilde
\eta_{j}<2\varepsilon_{n-1}^{3/4}=\eta.
\end{eqnarray*}
Moreover, for $\nu\leq N$,
we have
\begin{eqnarray*}
\lefteqn{K_\nu\leq K_N\leq \frac{1}{2\pi\sigma_{N}}\ln\frac{1}{\tilde
\eta_{N-1}}+1  =
\frac{4\cdot 3^{N-1}Q_n^3}{\pi r_0}\ln\frac{1}{2\varepsilon_{n-1}}+1}\\
&<& \frac{2Q_n^{\frac{2^nc
U}{24}+3}}{r_0}(\frac{1}{\varepsilon_{n-1}})^{\frac{1}{12\tau_*}}
<\left(\frac{\gamma'\sigma_N}{ 32\pi\cdot
2\varepsilon_{n-1}^{3/4}}\right)^{\frac{1}{\tau'}}.
\end{eqnarray*}
Then by the assumption $\rho_f\in D_{\alpha}(\gamma',\tau')$,
we can apply Proposition \ref{homo-lemma} to (\ref{homo-matri-2}),
getting an approximate solution $h_\nu^1$, with the error term $\tilde
P_\nu^{1}$, satisfying the estimations
\begin{equation}\label{est-erro}\begin{array}{ll}
&\|h_\nu^1\|_{\tilde
r_{\nu-1}-\sigma_{\nu},\tilde s_{\nu-1}-\delta_{\nu}}\leq \frac{C_0(d)K_\nu^{\tau'}\tilde
\eta_{\nu-1}}{\gamma'\sigma_\nu^{d}}<\tilde
\eta_{\nu-1}^{4/5},\\
&\|\tilde P_\nu^1\|_{ \tilde
r_\nu,\tilde s_\nu}<\frac{C_1(d)K_\nu^{d}\tilde
\eta_{\nu-1}^2}{\sigma_\nu^{d-1}}<\tilde
\eta_{\nu-1}^{7/4}.\end{array}
\end{equation}
Denote
$$h_\nu^\iota(\varphi,\theta)=\sum_{k\in \Z, l\in \Z^{d-1}} \hat {h}_{\nu,l}^{\iota}(k) e^{ 2\pi i k\varphi+ 2\pi i\la l,\theta\ra}, \ \ \iota=1,2.$$
If we let $\hat{h}_{\nu,l}^2(k)=\overline{\hat{h}_{\nu,-l}^1(-k)}$, then by the
relation between $w_{\nu-1}^1$ and $w_{\nu-1}^2$,  as well as (\ref{homo-matri-2}) and
(\ref{homo-matri-3}), we obtain that $h_\nu^2$ is an approximate
solution of (\ref{homo-matri-3}), with the error term $\tilde
P_\nu^{2}$, satisfying
the same estimates as in (\ref{est-erro}). Moreover, it holds that $h_\nu^1+h_\nu^2, i(h_\nu^1-h_\nu^2),
 \tilde P_\nu^1+
\tilde P_\nu^2, i(\tilde P_\nu^1-\tilde P_\nu^2)\in C_{\tilde r_\nu,\tilde s_\nu}^\omega(\T^1\times \T^{d-1},\R)$.
Then we have
\begin{eqnarray*}
&&A_{\nu-1}^{-1}Y_{\nu}(\cdot +\alpha)A_{\nu-1}- Y_{\nu}+\tilde F_{\nu-1}\\
&=& \left(
\begin{array}{ccc}
   0 &   [\tilde F_{\nu-1}^-]_\theta \cr  -[\tilde F_{\nu-1}^-]_\theta & 0\end{array} \right)+\left(
\begin{array}{ccc}
   \frac{1}{2}(\tilde P_\nu^1+\tilde P_\nu^2) &  \frac{i}{2}(\tilde P_\nu^1-\tilde P_\nu^2) \cr
   \frac{i}{2}(\tilde P_\nu^1-\tilde P_\nu^2)
   & -\frac{1}{2}(\tilde P_\nu^1+\tilde P_\nu^2)\end{array}
   \right)\\
&=:& \left(
\begin{array}{ccc}
   0 &   [\tilde F_{\nu-1}^-]_\theta \cr  -[\tilde F_{\nu-1}^-]_\theta & 0\end{array}
   \right)+\tilde P_\nu.
\end{eqnarray*}
In conclusion, we obtain that $Y_\nu, \tilde P_\nu \in C_{\tilde r_\nu, \tilde
s_\nu}^\omega(\T^1\times \T^{d-1}, sl(2,\R))$, with
\[\|Y_\nu\|_{\tilde r_{\nu-1}-\sigma_\nu, \tilde s_{\nu-1}-\delta_\nu}<\tilde \eta_{\nu-1}^{4/5},
\quad
 \|\tilde P_\nu\|_{\tilde r_\nu,\tilde s_\nu}<\tilde
\eta_{\nu-1}^{7/4}.\]
Let
\[\tilde
A_\nu=\tilde A_{\nu-1}e^{[\tilde F_{\nu-1}^{-}]_\theta J }, \ \textrm{and}\  \ E_\nu=\tilde A_\nu^{-1}e^{Y_\nu(\cdot+\alpha)}\tilde A_{\nu-1}e^{F_{\nu-1}}e^{-Y_\nu}. \]
Then by standard estimation (c.f. \cite{E92} for example), we have
\begin{eqnarray*}
\|E_\nu-id\|_{\tilde r_\nu, \tilde s_\nu}&\leq& C (\|\tilde F_{\nu-1}\|^2
+ \|Y_\nu\|^2
+\|\tilde F_{\nu-1}\|\cdot \|Y_\nu\|+\|\tilde P_\nu\|)\\
&\leq& C\tilde \eta_{\nu-1}^{8/5}\ll 1.
\end{eqnarray*}
Thus, by the implicit function theorem, there exists $\tilde F_\nu\in C^\omega_{\tilde r_\nu, \tilde s_\nu}(\T^1\times \T^{d-1}, sl(2,\R))$ such that $E_\nu=e^{\tilde F_\nu}$ with $\|\tilde F_\nu\|_{\tilde r_\nu, \tilde s_\nu}\leq \|E_\nu-id\|_{\tilde r_\nu, \tilde s_\nu}<\tilde \eta_{\nu-1}^{3/2}$, which completes the $\nu-$th step.

%

To finish the proof, we need to give the estimate of $\tilde F_N$.
By the selection of $U$, we have $Q_n^U\geq \ln Q_{n+1}$.
Therefore,
\begin{eqnarray*}
\lefteqn{\|\tilde F_{N}\|_{\tilde{r}_{N},\tilde{s}_{N}}\leq
\tilde\eta_0^{(\frac{3}{2})^{N}}
=\tilde \eta_0 e^{-((\frac{3}{2})^{N}-1)\ln\frac{1}{\tilde\eta_0}}}\\
&\leq&\tilde\eta_0 e^{-Q_n^{2^{n-1}c_1\tau_* U \ln (3/2)}2^nc\tau_* U\ln 2}
<\frac{\tilde\eta_0}{4Q_{n+1}^{2^{n+1}c\tau_* U}}=\varepsilon_{n}/2.
\end{eqnarray*}
Furthermore, it is obvious that $\tilde s_N>s_n, \tilde r_N>r_n$.
Let $$y_+:=y_N, \ \ \tilde g_+:=\tilde g_N,\ \ \tilde F_+:=\tilde
F_N.$$ Then we get the system (\ref{sys-after-iter}) with $\|\tilde g_+-\tilde g\|_{r_n}<2\tilde\eta_0$. Now, we give the estimate of the conjugation map
$\Psi_n:=e^{Y_N}e^{Y_{N-1}}\cdots e^{Y_1}$, satisfying
\begin{eqnarray*}
\|\Psi_n-I\|_{r_n,s_n}\leq \sum_{j=1}^N \|(e^{Y_j}-I)e^{Y_{j-1}}\cdots
e^{Y_2}e^{Y_1}\|_{\tilde r_j,\tilde s_j}\leq \sum_{j=1}^N4\|Y_j\|_{\tilde r_j, \tilde
s_j}<\varepsilon_{n-1}^{3/4}.
\end{eqnarray*}
Moreover, by Lemma \ref{rota-number-rem}, since $\Psi_n$ is close to the identity, meaning that $\Psi_n$ is homotopic to the identity, we have the fibered rotation number of (\ref{sys-after-iter}) is same as (\ref{sys-after-reso}) and also (\ref{sys-origin}).
\end{pf}

 As mentioned in Remark \ref{large-remark}, the transformation
 $e^{v_nJ}$ may not be close to the identity, which will be an
 obstruction on the convergence of the transformations. Therefore,
 we will do one more operation, which guarantees the convergence of
 the transformations.
 \begin{Lemma}\label{end-step-lemma}
Under the assumptions of Lemma \ref{redu-lemma}, the cocycle
(\ref{sys-after-iter}) can be conjugate to
\begin{equation}\label{sys-new}
(\alpha,\,R_{\varrho+ \frac{g_+(\varphi)}{2\pi}}e^{ F_+(\varphi,\theta)})\in
\mathcal
F_{r_{n}, s_{n}}(4\tilde \varepsilon_{n},\varepsilon_{n},\gamma', \tau')
\end{equation}
by the conjugation $R_{-\frac{v_n(\varphi)}{2\pi}}$, where $v_n(\varphi)$ is the same as
in Lemma \ref{reson-lemma}.

 \end{Lemma}

\begin{pf}
By simple computations, we know that $e^{v_n(\varphi)J}$ conjugates the cocycle (\ref{sys-after-iter}) to $(\alpha,\,R_{ \rho_f+  \frac{\tilde g_++\mathcal T_{Q_n}g-\hat g(0)}{2\pi}}e^{v_n J}e^{\tilde F_+}e^{-v_nJ})$.
Let
\[ g_+(\varphi)=2\pi(\rho_f-\varrho)+\tilde g_++\mathcal T_{Q_n}g-\hat g(0),\ \ \ \ F_+(\varphi,\theta)=e^{v_n J}\tilde F_+e^{-v_nJ}. \]
Recalling that
\[  v_n(\varphi+\tilde \alpha)-v_n(\varphi)=
-\mathcal T_{Q_n} g(\varphi)+\hat g(0),\]
we thus have
\begin{eqnarray*}
g_+-g=2\pi (\rho_f-\varrho)+ \tilde g_++\mathcal T_{Q_n}g-\hat g(0)-g=\tilde g_+-\tilde g+\mathcal R_{Q_n}g+\mathcal T_{Q_n}g-g,
\end{eqnarray*}
which implies that $\|\tilde g-g\|_{r_n}\leq 4\varepsilon_{n-1}$ by Lemma \ref{redu-lemma}.
And Lemma \ref{reson-lemma} implies
\[\|F_+\|_{r_n, s_n}\leq e^{\|\Im v_n\|_{r_n}}\|\tilde F_+\|_{r_n, s_n}e^{\|\Im v_n\|_{r_n}}<\varepsilon_n.\]
Moreover, by the same reason of Lemma \ref{reson-lemma}, the fibered rotation number of (\ref{sys-new}) does not change.

\end{pf}

Now we are on the position to prove Proposition \ref{iter-lemma}.\medskip

\textit{Proof of Proposition \ref{iter-lemma}.} Let
$\Phi_n=e^{v_nJ}\Psi_ne^{-v_nJ}$. Combining Lemma
\ref{reson-lemma}, \ref{redu-lemma}, and \ref{end-step-lemma}, we only need to
estimate $\Phi_n-I$:
\begin{eqnarray*}
\lefteqn{\|\Phi_n-I\|_{r_n,s_n}=\|e^{v_nJ}\Psi_ne^{-v_nJ}-I\|_{r_n,s_n}}\\
&=&\|e^{-v_nJ}(\Psi_n-I)e^{v_nJ}\|_{r_n,s_n}\\&\leq& 2\|\Psi_n-I\|_{r_n,s_n}
<2\varepsilon_{n-1}^{3/4}.
\end{eqnarray*}
\qed

\subsection{Proof of Theorem \ref{thm-ana-tech}}

Let $\varepsilon_0=\|F\|_{\mathfrak r}$ satisfy (\ref{pert-condi}), where $r_0=s_0=\mathfrak r$,
and $(Q_n)$ is the selected sequence of $\tilde\alpha$ in Lemma
\ref{cd-coroll} with $\mathcal A=4$, $\tilde U=\tilde U(\tilde\alpha)$ and $U:=\tilde U+8$.
Define $r_n, s_n, \varepsilon_n, \tilde\varepsilon_n$ for $n\geq 1$ as in (\ref{para-induc}). Let $\rho_f:=\rho(\alpha,R_{\varrho_0}e^F)$, $g_0:=2\pi (\varrho_0-\rho_f)$ and $F_0:=F$. By the assumption,
we have $(\alpha, R_{\rho_f+\frac{g_0}{2\pi} }e^{F_0})\in \mathcal{F}_{r_0,s_0}(8\varepsilon_0,\varepsilon_0,\gamma',\tau')$. Then by Lemma \ref{redu-lemma}, there exists $\Psi_1\in C^\omega_{r_1, s_1}(\T^1\times\T^{d-1}, SL(2,\R))$ that conjugates $(\alpha, R_{\rho_f+\frac{g_0}{2\pi}}e^{F_0})$ to $(\alpha, R_{\rho_f+\frac{g_1}{2\pi}}e^{F_1})$ which belongs to $\mathcal F_{r_1, s_1}(4\tilde\varepsilon_1, \varepsilon_1, \gamma', \tau')$. Let $\Phi_1:=\Psi_1$.
Applying Proposition \ref{iter-lemma} inductively for $n\geq 2$, then we obtain a sequence
of transformations $\{\Phi_n\}_{n\geq 2}$ such that $\Phi_n$
conjugates (\ref{sys-iter-former}) into (\ref{sys-iter-after}) with $\varrho=\rho_f$.
Let \[\Phi^{(n)}=\Phi_n\Phi_{n-1}\cdots \Phi_1.\]
Then by the estimates in Proposition \ref{iter-lemma}, we get
\begin{eqnarray*}
\|\Phi^{(n+1)}-\Phi^{(n)}\|_{
r_{n+1},s_{n+1}}&=&\|(\Phi_{n+1}-I)\Phi^{(n)}\|_{r_{n+1},s_{n+1}} \\
&\leq&
2\varepsilon_n^{3/4}\prod_{j=1}^{n-1}(1+2\varepsilon_j^{3/4})<4\varepsilon_n^{3/4}.
\end{eqnarray*}
Then $\Phi^{(n)}$ and $g_n$ converge to $\Phi, g_\infty$ respectively in $C^0$ topology. Hence,
$\Phi$ conjugates
$(\alpha, R_{\varrho_0}e^F)$ to the cocycle $(\alpha,\, R_{\rho_f+\frac{g_\infty(\varphi)}{2\pi}})$.
Now, we certificate that $\Phi(\cdot)$ is actually in $C^\infty$, that is $\{\Phi^{(n)}(\cdot)\}_{n\in\N}$ converges
in $C^\infty$ topology.
By the definition of $\{\varepsilon_n\}_{n\in \N}$, for any
$j\in \N^{d}$, there exists $n_j\in\N$ such that for any $n\geq n_j$, we
have $Q_n^{-1}< r_0/2$ and  $Q_n^{4|j|}<\varepsilon_{n-1}^{-1/2}$,
which implies
$$(2r_0^{-1}Q_n^3)^{|j|}\varepsilon_{n-1}^{3/4}<\varepsilon_{n-1}^{1/4}, \ \forall n\geq n_j.$$
Then for any $j\in
\N^{d}, n\geq n_j$, by Cauchy estimates, we have
$$\left|D^j(\Phi^{(n+1)}-\Phi^{(n)})\right|\leq \frac{\|\Phi^{(n+1)}-\Phi^{(n)}\|_{r_{n+1},s_{n+1}}}
{r_{n+1}^{|j|}}\leq \left(\frac{2Q_{n+1}^3}{r_0}\right)^{|j|}\cdot
4\varepsilon_{n}^{3/4}<4\varepsilon_{n}^{1/4}.$$ Hence, the limit
$\Phi$ is in $C^\infty$, which also means $g_\infty\in
C^\infty(\T^1,\R)$.
\qed

\section{Proof of Theorem \ref{thm-2-fre-tech}}\label{sec:6}

In this section, we also use KAM iterations to prove   Theorem \ref{thm-2-fre-tech}. Similar as in the proof of Theorem \ref{thm-ana-tech}, the key point in KAM theory is to solve the linearized homological equation, which involves the small divisor problems. Therefore, we first give some important small divisor lemmas in our proof, then perform the KAM scheme, and finally verify the convergence of the iterations.

\subsection{Small divisor lemmas}
Let $\alpha=(\tilde\alpha,\alpha')\in\T^1\times \T^{1}$. In the sequel, we will fix some $\chi\geq 5$ and let $(\tilde\alpha,\alpha')\in\Omega(\chi)$. Assume that $\frac{\tilde p_n}{\tilde q_n}$ and $\frac{p_n'}{q_n'}$ are  the continued fraction approximates to $\tilde \alpha$ and $\alpha'$
respectively. 

The following two small divisor lemmas are important in our
proof.

\begin{Lemma}\label{lemma-2-fre-divisor}
Let $\chi\geq 5$. If $(\tilde\alpha, \alpha')\in \Omega(\chi)$, then for $0<|k|+|l|< \tilde q_n$ $(n\geq 1)$, we have
$$\|k\tilde\alpha+l \alpha'\|_{\R/\Z}\geq \frac{1}{2 q_n'}.$$
\end{Lemma}
\begin{pf}
We first deal with the case $kl=0$. If $l=0$, then for $0<|k|< \tilde q_{n}
$, we have
$$\|k\tilde\alpha\|_{\R/\Z}\geq \frac{1}{2\tilde q_{n}}>\frac{1}{2q_n'}.$$
If $k=0$, then for $0<|l|<\tilde q_n<q_n'$, we have
$$\|l \alpha'\|_{\R/\Z}\geq \frac{1}{2q_n'}.$$
Now we suppose $kl\neq 0$. Since $\tilde q_n\wedge q_{n-1}'=1, \tilde q_n\wedge
\tilde p_n=1, |k|<\tilde q_n$, then $\tilde q_n\nmid k \tilde p_n q_{n-1}'$, which implies that
$$\|l\frac{p_{n-1}'}{q_{n-1}'}+k\frac{\tilde p_n}{\tilde q_n}\|_{\R/\Z}\geq \frac{1}{q_{n-1}'\tilde q_n}.$$
Moreover,
$$|l(\alpha'-\frac{p_{n-1}'}{q_{n-1}'})|<\frac{\tilde q_n}{q_{n-1}'q_{n}'}<\frac{1}{4q_{n-1}'\tilde q_n},$$
$$|k(\tilde \alpha-\frac{\tilde p_n}{\tilde q_n})|<\frac{1}{\tilde q_{n+1}}<\frac{1}{4q_{n-1}'\tilde q_n}.$$
Therefore, we have
$$\|k\tilde \alpha+l\alpha'\|_{\R/\Z}\geq\frac{1}{2q_{n-1}'\tilde q_n}>\frac{1}{2q_n'}.$$
\end{pf}

One of the main observation for $(\tilde \alpha,\alpha')\in \Omega(\chi)$ is that
 if $\varrho\in D_\alpha(\gamma,\tau)$ for some $\gamma>0, \tau>2$, then the small divisors $\|\la k,\alpha\ra\pm 2\varrho\|_{\R/\Z}$ have some uniform lower bound for $|k|$ not too large, which is the content of the following lemma:
\begin{Lemma}\label{lemma-uniform-divisor}
Let $0<\gamma<1, \tau>2, \alpha=(\tilde \alpha,\alpha')\in\Omega(\chi)$, and  $\varrho\in D_\alpha(\gamma,\tau)$. Then for any $k\in\Z^2$, $|k|\leq \tilde q_{n+1}^{1/2} $, we have
\begin{equation}\label{beta-diophan}
\|\la k,\alpha\ra\pm 2\varrho\|_{\R/\
Z}\geq \frac{c\gamma^{\tau+1}}{q_n'^{\tau^2}},
\end{equation}
where $c$ is a constant only depending on $\tau$.
\end{Lemma}
\begin{pf}
Since $\|\la k,\alpha\ra-2\varrho\|_{\R/\Z}=\|\la -k, \alpha\ra+2\varrho\|_{\R/\Z}$, we only need to prove
\[\|\la k,\alpha\ra+2\varrho\|_{\R/\Z}\geq \frac{c\gamma^{\tau+1}}{q_n'^{\tau^2}}.\]
Let $L_n=\frac{12\cdot  (2 q_n')^\tau}{\gamma}$. We divide the proof into two cases:

\textit{Case 1: $\tilde q_{n+1}\leq L_n^2$.}  Then for $|k|\leq \tilde q_{n+1}^{1/2}\leq L_n$,  (\ref{beta-diophan}) follows directly from the  Diophantine condition that $\varrho\in D_\alpha(\gamma,\tau)$.

\textit{Case 2: $\tilde q_{n+1}>L_n^2$.}  Assume that
\[\|\la k,\alpha\ra+2\varrho\|_{\R/\Z}=|\la k,\alpha \ra+2\varrho+j_k|<1\]
 for some $j_k\in\Z$. Then we have $|j_k|<1+|\la k,\alpha\ra|+2|\varrho|\leq 3+|k|$. Let $k'=(j_k, k)\in\Z^3$, $\omega=(1,\alpha)$. Then
 \[|k'|=|j_k|+|k|<3+2|k|\leq 3\tilde q_{n+1}^{1/2}<\frac{3\tilde q_{n+1}}{L_n},\]
 because $\tilde q_{n+1}>L_n^2>144\cdot q_n'^{2\tau}$.
There exist $l,m\in\Z, \tilde k\in \Z^3$, such that
\[k'=l(-\tilde p_n,\tilde q_n,0)+m(- p_n', 0,q_n')+(\tilde k_1, \tilde k_2, \tilde k_3),\]
with $|\tilde k_2|<\tilde q_n, |\tilde k_3|<q_n', |l|\leq \frac{3\tilde q_{n+1}}{L_n\tilde q_n}, |m|\leq \frac{ 3\tilde q_{n+1}}{L_n  q_n'}$. Then,
\begin{eqnarray*}
\lefteqn{\|\la k,\alpha\ra+2\varrho\|_{\R/\Z}=|\la k', \omega\ra +2\varrho|}\\
&&=|l(\tilde\alpha \tilde q_n-\tilde p_n)+m( \alpha' q_n'- p_n')+\tilde k_1+\tilde k_2\tilde\alpha+\tilde k_3\alpha'+2\varrho|.
\end{eqnarray*}
Since $|\tilde k_2|+|\tilde k_3|<2q_n'$, together with the fact that $\varrho\in D_\alpha(\gamma,\tau)$, we get
\[\|\la\tilde k,\omega\ra+2\varrho\|_{\R/\Z}=\|\tilde k_2\tilde\alpha+\tilde k_3\alpha'+2\varrho\|_{\R/\Z}\geq \frac{\gamma}{(|\tilde k_2|+|\tilde k_3|+1)^\tau}\geq\frac{\gamma}{(2 q_n')^\tau}.\]
Therefore, we have
\begin{eqnarray*}
\lefteqn{\|\la k,\alpha\ra+2\varrho\|_{\R/\Z}=|\la\tilde k,\omega\ra+2\varrho+l(\tilde\alpha\tilde q_n-\tilde p_n)+m(\alpha'q_n'-p_n')|}\\
&\geq& \|\tilde k,\omega\ra+2\varrho\|_{\R/\Z}-|l(\tilde \alpha \tilde q_n-\tilde p_n)|-|m(\alpha'q_n'-p_n')|\\
&\geq& \frac{\gamma}{(2q_n')^\tau}-\frac{3\tilde q_{n+1}}{L_n \tilde q_n}\cdot \frac{1}{\tilde q_{n+1}}-\frac{3\tilde q_{n+1}}{L_n q_n'}\cdot\frac{1}{ q_{n+1}'}\geq\frac{\gamma}{2(2 q_n')^\tau}.
\end{eqnarray*}
This completes the proof.
\end{pf}\medskip

\subsection{KAM scheme}

For $\mathfrak r,\eta,\tilde \eta,\gamma>0,\tau>2$, we let
\begin{eqnarray*}
\lefteqn{\tilde \cF_{\mathfrak r}( \eta,\tilde \eta,\gamma,\tau) :=}\\
&&\left\{\begin{array}{ll} (\alpha,\,R_{\varrho+\frac{g}{2\pi} } e^{F})\ : g\in C_{\mathfrak r}^\omega(\T^2, \R), \ F\in C_{\mathfrak r}^\omega(\T^2, sl(2,\R)),
\\
 \rho(\alpha, R_{\varrho+\frac{g}{2\pi}})\in D_\alpha(\gamma,\tau),\ \textrm{and}\  \|g\|_{\mathfrak r}\leq \eta, \|F\|_{\mathfrak r}\leq \tilde \eta
  \end{array}
  \right\}.
  \end{eqnarray*}

Without loss of generality, we assume $\mathfrak r\leq 1$. For given $\chi\geq 5, 0<\mathfrak r_0:=\mathfrak r\leq 1,\tau>2$, we let $\ell_*\geq 2$ be the smallest real number such that for any $\ell\geq \ell_*$, we have
\begin{equation}
2^{11}\ell^{16\chi\tau^2}\leq e^{\ell^{1/10}\mathfrak r_0}.
\end{equation}
Let $\gamma>0$ be sufficiently small. Suppose $\varepsilon_0$ is small enough such that
\begin{equation}\label{condition-varepsilon}
\varepsilon_0<\min\{c\gamma^{8(\tau+1)}\mathfrak r_0^{16},\ \ell_*^{-16\chi\tau^2}\},
\end{equation}
where $c>0$ is a sufficiently small constant. Let $n_*\in\N$ be the smallest natural number such that for $n>n_*$
\begin{equation}\label{para-selec}
 2^{11}\tilde q_n^{16\chi\tau^2}\leq e^{\tilde q_n^{\frac{1}{10}}\mathfrak r_0}, \ \textrm{and}\ e^{-\tilde q_n^{\frac{1}{10}}\mathfrak r_0}\leq \varepsilon_0.
\end{equation}
Let $\bar {\mathfrak r}_0=\mathfrak r_0$. Now we inductively define the sequences for $n\geq 1$:
\begin{eqnarray}\label{para-induc-2-fre}
&& \mathfrak r_n=\frac{\mathfrak r_0}{2q_{n_*+n-1}'^2},\  \ \bar {\mathfrak r}_{n}=\frac{4\mathfrak r_0}{q_{n_*+n}'^2},\\ \nonumber
&& \varepsilon_n=\varepsilon_{n-1}e^{-\tilde q_{n_*+n}^{\frac{1}{10}}\mathfrak r_0/2},\ \ \ \ \tilde \varepsilon_n=4\sum_{j=0}^{n-1}\varepsilon_j^{3/4}.
\end{eqnarray}
Then by induction, we have
\begin{equation}\label{est-varepsilon}
e^{-\tilde q_{n_*+n+1}^{\frac{1}{10}}\mathfrak r_0}\leq \varepsilon_n.
\end{equation}

\medskip

\begin{Proposition}[Iterative Lemma]\label{prop-2-freq}
Let $\varepsilon_0, \mathfrak r_0, \gamma>0, \tau>2, \chi\geq 5$ satisfy (\ref{condition-varepsilon}),
 and $\alpha=(\tilde\alpha,\alpha')\in \Omega(\chi)$.  The following holds for $n\geq 1$:  If the  cocycle
\begin{equation}\label{sys-2-fre-before}
(\alpha,\, R_{\varrho+\frac{g_{n}}{2\pi}}e^{F_n})\in \tilde \cF_{\mathfrak r_{n}}( \tilde \varepsilon_{n}, \varepsilon_{n},\gamma,\tau)
\end{equation}
 satisfies $\mathcal R_{\tilde q_{n_*+n}} g_n=0$, then there exists $\Phi_n\in C_{\mathfrak r_{n+1}}^\omega(\T^2, SL(2,\R))$ with $\|\Phi_n- I\|_{\mathfrak r_{n+1}}<\varepsilon_{n}^{1/3}$ such that $x_{n+1}=\Phi_n x_{n}$ conjugates
 the cocycle (\ref{sys-2-fre-before}) to the cocycle
 \begin{equation}\label{sys-2-fre-after}
 (\alpha,\, R_{\varrho+\frac{g_{n+1}}{2\pi}}e^{F_{n+1}})\in \tilde {\mathcal{F}}_{\mathfrak r_{n+1}}(\tilde \varepsilon_{n+1},\varepsilon_{n+1},\gamma,\tau)
 \end{equation}
 with $\|g_{n+1}-g_n\|_{\mathfrak r_{n+1}}\leq4\varepsilon_n^{3/4}$, $\|F_{n+1}\|_{\mathfrak r_{n+1}}\leq \varepsilon_{n+1}$, and
 $\mathcal R_{\tilde q_{n_*+n+1}}g_{n+1}=0$.
\end{Proposition}

We divide the proof into three parts that corresponds to Lemma \ref{lemma-2-fre-nonre}, \ref{lemma-2-fre-redu} and \ref{lemma-2-fre-revers},
and we omit the subscript for convenience in case there is no confusing. \medskip

\begin{Lemma}\label{lemma-2-fre-nonre}
For $n\geq 1$, if $(\alpha,\,R_{\varrho+\frac{g}{2\pi}}e^{F})\in\tilde {\mathcal{F}}_{\mathfrak r_n}(\tilde \varepsilon_n,\varepsilon_n,\gamma,\tau)$ with $\mathcal R_{\tilde q_{n_*+n}}g=0$,
then  it can be conjugate to $(\alpha, R_{\rho_f }e^{\tilde F})\in\tilde{\mathcal{ F}}_{\bar {\mathfrak r}_n}( 0,60\varepsilon_n,\gamma,\tau)$, with $\rho_f=\rho(\alpha, R_{\varrho+\frac{g}{2\pi}}e^F)$.
\end{Lemma}

\begin{pf}
Let $v$ be the solution of
\begin{equation}\label{equ-2-fre-v}
 v( \phi+ \alpha)-v(\phi)=-g(\phi)+\hat g(0).
\end{equation}
Then $e^{-v(\phi)J}$ conjugates  $(\alpha,\,R_{\varrho+\frac{g}{2\pi}}e^{F})$ to
\begin{equation}
(\alpha,\, R_{\varrho+\frac{\hat g(0)}{2\pi}}e^{e^{-vJ} F( \phi)e^{vJ}}).
\end{equation}
Applying Lemma \ref{claim-radius}, we get
\[\|\mathfrak{I}v\|_{\bar {\mathfrak r}_n}\leq \frac{Cq_{n_*+n-1}'^6\varepsilon_0^{3/4} }{{\mathfrak r}_0^2q_{n_*+n}'}<\varepsilon_0^{1/2}\ll 1,\]
 because  $\|\la k,\alpha\ra\|_{\R/\Z}\geq \frac{1}{2q_{n_*+n}'}$ for $0<|k|<\tilde q_{n_*+n}$ by Lemma \ref{lemma-2-fre-divisor}.  Then we have
\[\|e^{-vJ} Fe^{vJ}\|_{\bar {\mathfrak r}_n}\leq \|e^{vJ}\|_{\bar {\mathfrak r}_n}\|F\|_{\bar {\mathfrak r}_n}\|e^{-vJ}\|_{\bar {\mathfrak r}_n}\leq e^{2\|\mathfrak{I}v\|_{\bar {\mathfrak r}_n}}\|F\|_{{\mathfrak r}_n}<2\varepsilon_n.\]
By Lemma \ref{Lem-rot-num},  we have
\begin{equation}\label{equ-rota-diff}
|\varrho+\frac{\hat g(0)}{2\pi}-\rho_f|\leq \|e^{e^{-vJ}Fe^{vJ}}-I\|_0\leq 2\|e^{-vJ}Fe^{vJ}\|_{0}<4\varepsilon_n.
\end{equation}
Let $E=R_{\varrho+\frac{\hat g(0)}{2\pi}-\rho_f}e^{e^{-vJ}Fe^{vJ}}$. Then
\begin{eqnarray*}
\lefteqn{\|E-I\|_{\bar {\mathfrak r}_n}\leq \|e^{(2\pi(\varrho-\rho_f)+\hat g(0))J}-I\|+\|e^{e^{-vJ}Fe^{vJ}}-I\|_{\bar {\mathfrak r}_n}}\\
&\leq & 4\pi|\varrho+\frac{\hat g(0)}{2\pi}-\rho_f|+2\|e^{-vJ}Fe^{vJ}\|_{\bar {\mathfrak r}_n}<60\varepsilon_n.
\end{eqnarray*}
By the implicit function theorem, there exists $\tilde F\in C_{\bar {\mathfrak r}_n}^\omega(\T^2, sl(2,\R))$ such that $E=e^{\tilde F}$ with $\|\tilde F\|_{\bar {\mathfrak r}_n}\leq \|E-I\|_{\bar {\mathfrak r}_n}<60\varepsilon_n$.
Moreover, since the conjugation map $R_{\frac{v}{2\pi}}$ is homotpoic to the identity, the fibered rotation number of the new cocycle stays the same by Lemma \ref{rota-number-rem}.
\end{pf}

Denote by $\mathcal B_{\mathfrak r}^\omega(\T^2, sl(2,\R))$ the set of $F\in C^\omega(\T^2, sl(2,\R))$  satisfying
\[|F|_{\mathfrak r}:=\sum_{k\in\Z^2}\|\hat F(k)\|e^{2\pi|k|\mathfrak r}<\infty.\]
Then for any $F\in C_{\mathfrak r}^\omega(\T^2, sl(2,\R))$, we have $F\in \mathcal B_{{\mathfrak r}_+}$ for any $0<{\mathfrak r}_+<\mathfrak r$ with the estimate
\[|F|_{{\mathfrak r}_+}\leq \frac{36}{\min\{1, (\mathfrak r-{\mathfrak r}_+)^2\}}\|F\|_{\mathfrak r}.\]

Now, for $\mathfrak r>0, A\in SL(2,\R), \eta>0$, we  decompose 
\[\mathcal B_{\mathfrak r}= \mathcal B_{ \mathfrak r}^{(nre)}(\eta)\oplus \mathcal B_{ \mathfrak r}^{(re)}(\eta),\]
 such that for any $Y\in \mathcal B_{ \mathfrak r}^{(nre)}(\eta)$,
\begin{equation}\label{condition-nonresonant}
A^{-1}Y(\cdot+\alpha)A-Y\in \mathcal B_{ \mathfrak r}^{(nre)},\ \  |A^{-1}Y(\cdot+\alpha)A-Y|_{ \mathfrak r}\geq \eta|Y|_{\mathfrak r}.
\end{equation}

\begin{Lemma}\label{lemma-2-fre-redu}
For $n\geq 0$, if $(\alpha,\, R_{\rho_f}e^{\tilde F})\in\tilde {\mathcal F}_{\bar {\mathfrak r}_n}(0, 60\varepsilon_n,\gamma,\tau)$,
then there exists a map $\Psi\in C_{\bar {\mathfrak r}_n/2}^\omega(\T^2, SL(2,\R))$ with $\|\Psi-I\|_{\bar {\mathfrak r}_n/2}<4\varepsilon_n^{3/8}$, such that
the cocycle $(\alpha, R_{\rho_f}e^{\tilde F})$ can be conjugated  to
\begin{equation}\label{sys-2-fre-after-3}
(\alpha,\, R_{\rho_f+\frac{\bar g(\phi)}{2\pi}}e^{\bar F(\phi)})\in \tilde {\mathcal{F}}_{\bar {\mathfrak r}_n/4}(120\varepsilon_n,\varepsilon_{n+1}/2,\gamma,\tau),
\end{equation}
with $\mathcal R_{\tilde q_{n_*+n+1}}\bar g=0$, where $\rho_f=\rho(\alpha, R_{\rho_f}e^{\tilde F})$.
\end{Lemma}

In order to prove Lemma \ref{lemma-2-fre-redu}, we will need the following lemma \medskip

\begin{Lemma}[{\cite{CCYZ}, Lemma 3.1}]\label{lemma-houyou}
Let $\varepsilon\leq (4\|A\|)^{-4}$ and $\eta\geq13\|A\|^2\varepsilon^{1/2}$. Then for any $F\in \mathcal B_{\mathfrak r}$ with $|F|_{\mathfrak r}\leq \varepsilon$, there exist $Y\in \mathcal B_{\mathfrak r}, F^{(re)}\in \mathcal B_{\mathfrak r}^{(re)}(\eta)$ such that
\[   e^{Y(\cdot+\alpha)} Ae^F e^{-Y}=Ae^{F^{(re)}},\]
with $|Y|_{\mathfrak r}\leq \varepsilon^{1/2}$ and $|F^{(re)}|_{\mathfrak r}\leq 2\varepsilon$.
\end{Lemma}\medskip

\begin{pf}(Proof of Lemma \ref{lemma-2-fre-redu}).

For $n=0$ we have $|\tilde F|_{\bar {\mathfrak r}_0/2}\leq C{\mathfrak r}_0^{-2}\varepsilon_{0}<\varepsilon_0^{3/4}$, and for $n\geq 1$,
\[|\tilde F|_{\bar {\mathfrak r}_n/2}\leq C{\mathfrak r}_0^{-2}\tilde q_{n_*+n}^{4\chi}\varepsilon_{n}\leq \varepsilon_n^{3/4}\]
by (\ref{para-selec}). Denote $\epsilon_n:=\varepsilon_n^{3/4}$.
Assume $\Lambda$ is a subset of $\Z^2$ with $\Lambda=-\Lambda$ such that
\[k\in\Lambda\Rightarrow \|\la k,\alpha\ra \pm 2\rho_f \|_{\R/\Z}\geq 13\epsilon_n^{1/2}.\]

Define $\mathcal B_{\bar {\mathfrak r}_n/2}^{1}$ as the space of all $ F\in \mathcal B_{\bar {\mathfrak r}_n/2}$ of the form
\[ F( \phi)=\sum_{k\in\Lambda}M^{-1}\left(\begin{array}{cc} 0 & \hat{F}^{11}(k)-i\hat{ F}^+(k)\cr
\hat{F}^{11}(k)+i\hat{ F}^+(k) & 0\end{array}\right)Me^{2\pi i\la k, \phi\ra}\]
and $\mathcal B_{\bar {\mathfrak r}_n/2}^{2}$ as the space of all $ F\in\mathcal B_{\bar {\mathfrak r}_n/2}$ of the form
\begin{eqnarray}\label{equ-reson-form}
&&F(\phi)=\sum_{k\in\Z^2}M^{-1}\left(\begin{array}{cc} i\hat F^-(k) & 0 \cr
0 & -i\hat F^-(k)\end{array}\right)Me^{2\pi i\la k,\phi\ra}\\
&+&\sum_{k\in\Z^2\backslash\Lambda}M^{-1}\left(\begin{array}{cc} 0 & \hat{ F}^{11}(k)-i\hat{ F}^+(k)\cr
\hat{ F}^{11}(k)+i\hat{ F}^+(k) & 0\end{array}\right)Me^{2\pi i\la k,\phi\ra}\nonumber,
\end{eqnarray}
where $M=\frac{1}{1+i}\left(\begin{array}{cc} 1 & -i\cr
1 & i\end{array}\right)$ and $\hat{ F}^\pm(k)=\frac{1}{2}(\hat{ F}^{12}(k)\pm \hat{ F}^{21}(k))$.
It is obvious that $\mathcal B_{\bar {\mathfrak r}_n/2}=\mathcal B_{\bar {\mathfrak r}_n/2}^{1}\oplus\mathcal B_{\bar {\mathfrak r}_n/2}^{2}$.
Now, we denote $A=R_{\rho_f}$.
One can check that for $Y\in \mathcal B_{\bar {\mathfrak r}_n/2}^1$,
\[|A^{-1}Y(\cdot+\alpha)A-Y|_{\bar {\mathfrak r}_n/2}\geq 13\epsilon_n^{1/2}|Y|_{\bar {\mathfrak r}_n/2},\]
and thus $\mathcal B_{\bar {\mathfrak r}_n/2}^1\subseteq\mathcal B_{\bar {\mathfrak r}_n/2}^{(nre)}(13\epsilon_n^{1/2})$.
Then, we can apply Lemma \ref{lemma-houyou}, getting some $\tilde Y\in\mathcal B_{\bar {\mathfrak r}_n}$ and $\tilde F^{(re)}$ in the form
(\ref{equ-reson-form}) such that $e^{\tilde Y}$ conjugates $(\alpha,\, Ae^{\tilde F})$ to $(\alpha,\,Ae^{\tilde F^{(re)}})$
with $\|\tilde Y\|_{\bar {\mathfrak r}_n/2}\leq \epsilon_n^{1/2}$ and $\|\tilde F^{(re)}\|_{\bar {\mathfrak r}_n/2}\leq 2\epsilon_n$.
 Moreover, since $\rho_f\in D_\alpha(\gamma,\tau)$, then for any $|k|\leq \tilde q_{n_*+n+1}^{1/2}$ by Lemma \ref{lemma-uniform-divisor}, we have
\[\|\la k,\alpha\ra\pm 2\rho_f\|_{\R/\Z}\geq \frac{c\gamma^{\tau+1}}{q_{n_*+n}'^{\tau^2}}>13\epsilon_n^{1/2}.\]
The last inequality obviously holds because of (\ref{condition-varepsilon}), (\ref{para-selec}) and the selection of $\ell_*$. It implies that $\Lambda^c\subseteq \{k\in\Z^2\ |\ |k|>\tilde q_{n_*+n+1}^{1/2}\}$.
We now let 
\[\bar g(\phi) =-\sum_{|k|<\tilde q_{n_*+n+1}}\hat{\tilde F}^{(re)-}(k)e^{2\pi i\la k, \phi\ra}, \ \ \textrm{and}\ E=e^{\bar g J}e^{\tilde F^{(re)}}=R_{-\frac{\bar g}{2\pi}}e^{\tilde F^{(re)}}.\]
Then we have $\mathcal T_{q_{n_*+n+1}^{1/2}}(\tilde F^{(re)}+\bar gJ)=0$, which implies
 \begin{eqnarray*}
 \|\tilde F^{(re)}+\bar g J\|_{\bar {\mathfrak r}_n/4}&\leq& C\tilde q_{n_*+n+1}e^{-\tilde q_{n_*+n+1}^{1/2}\bar {\mathfrak r}_n/4}\|\tilde F^{(re)}+\bar g J\|_{\bar {\mathfrak r}_n}\\
 &<&\frac{\varepsilon_ne^{-\tilde q_{n_*+n+1}^{1/2}\bar {\mathfrak r}_n/8}}{4}\leq\varepsilon_{n+1}/4
 \end{eqnarray*}
 by (\ref{para-selec}) and (\ref{est-varepsilon}).
Then,
\begin{eqnarray*}
\|E-I\|_{\frac{\bar {\mathfrak r}_n}{4}}&=& \|e^{\bar g J}e^{-\bar g J+\tilde F^{(re)}+\bar gJ}-e^{\bar g J}e^{-\bar gJ}\|_{\frac{\bar {\mathfrak r}_n}{4}}\\
&\leq& \|e^{\bar g J}\|_{\frac{\bar {\mathfrak r}_n}{4}}\cdot \|e^{-\bar gJ+\tilde F^{(re)}+\bar gJ}- e^{-\bar gJ}\|_{\frac{\bar {\mathfrak r}_n}{4}}\\
&\leq & e^{2\|\bar g J\|_{\frac{\bar {\mathfrak r}_n}{4}}}e^{\|\tilde F^{(re)}+\bar gJ\|_{\frac{\bar {\mathfrak r}_n}{4}}}\|\tilde F^{(re)}+\bar gJ\|_{\frac{\bar {\mathfrak r}_n}{4}} \\
&\leq & 2\|\tilde F^{(re)}+\bar g J\|_{\bar {\mathfrak r}_n/4}\leq \varepsilon_{n+1}/2.
\end{eqnarray*}
By the implicit function theorem, there exists $\bar F\in C_{\bar {\mathfrak r}_n/2}^\omega(\T^2, sl(2,\R))$ such that $E=e^{\bar F}$ with $\|\bar F\|_{\bar {\mathfrak r}_n/4}\leq \|E-I\|_{\bar {\mathfrak r}_n/4}\leq \varepsilon_{n+1}/2$.

Let $ \Psi=e^{\tilde Y}$. Then $\|\Psi-I\|_{\bar {\mathfrak r}_n/2}\leq 2\|\tilde Y\|_{\bar {\mathfrak r}_n/2}\leq 2\varepsilon_n^{3/8}$, and $ \Psi$
conjugates
the cocycle $(\alpha, R_{\rho_f}e^{\tilde F})$ to $(\alpha, R_{\rho_f+\frac{\bar g(\phi)}{2\pi}}e^{\bar F(\phi)})$, with $\|\bar g\|_{\bar {\mathfrak r}_n/4}\leq 120\varepsilon_n$.

To that end, since $\Psi$ is close to the identity and thus homotopic to the identity, Lemma \ref{rota-number-rem} implies that the fibered rotation number of the new cocycle is the same as the original one.
\end{pf}

\begin{Lemma}\label{lemma-2-fre-revers}
Under the assumptions of Lemma \ref{lemma-2-fre-redu}, the cocycle
$(\alpha, R_{\rho_f+\frac{\bar g}{2\pi}}e^{\bar F})$ can be conjugated to $(\alpha, R_{ \varrho+\frac{g_+}{2\pi}}e^{F_+})\in \mathcal {\tilde F}_{{\mathfrak r}_{n+1}}(\tilde\varepsilon_{n+1}, \varepsilon_{n+1},\gamma,\tau)$.
\end{Lemma}
\begin{pf}

 With the conjugation $e^{v (\phi)J}$, where $v$ is the same as in Lemma \ref{lemma-2-fre-nonre}, we get the cocycle $(\alpha, R_{\rho_f+\frac{\bar g+g-\hat g(0)}{2\pi}}e^{vJ}e^{\bar F}e^{-vJ})$.
Then similar to the proof in Lemma \ref{lemma-2-fre-nonre}, we let
\[g_+=2\pi(\rho_f-\varrho)+\bar g+g-\hat g(0),\ \ \textrm{and}\ F_+=e^{vJ}\bar Fe^{-vJ}.\]
By Lemma \ref{claim-radius},  Lemma \ref{lemma-2-fre-redu} and estimation (\ref{equ-rota-diff}), we get
\[\|g_+-g\|_{{\mathfrak r}_{n+1}}\leq 2\pi|\varrho+\frac{\hat g(0)}{2\pi}-\rho_f|+\|\bar g\|_{{\mathfrak r}_{n+1}}<\varepsilon_n^{3/4},\]
\[\|F_+\|_{{\mathfrak r}_{n+1}}\leq \|e^{vJ}\|_{{\mathfrak r}_{n+1}}^2\cdot \|\bar F\|_{{\mathfrak r}_{n+1}}<\varepsilon_{n+1}.\]
Moreover, the fibered rotation number does not change because the conjugation map $e^{vJ}$ is homotopic to the identity.
 \end{pf}

\textit{Proof of Proposition \ref{prop-2-freq}:} Let $\Phi_n= e^{v (\phi)J}\Psi e^{-v(\phi)J}$, where $v$ is defined as in (\ref{equ-2-fre-v}) and $\Psi$ is the map in Lemma \ref{lemma-2-fre-redu}. Then combining Lemma \ref{lemma-2-fre-nonre}, \ref{lemma-2-fre-redu}, \ref{lemma-2-fre-revers},  the cocyle (\ref{sys-2-fre-before}) is conjugate to (\ref{sys-2-fre-after}) via the conjugation $\Phi_n$.  Moreover,
\begin{eqnarray*}
\|\Phi-I\|_{{\mathfrak r}_{n+1}}=\|e^{vJ}(\Psi-I)e^{-vJ}\|_{{\mathfrak r}_{n+1}}\leq 2\|\Psi-I\|_{\bar {\mathfrak r}_n/2}\leq 8\varepsilon_n^{3/8}<\varepsilon_n^{1/3}.
\end{eqnarray*}

\subsection{Proof of Theorem \ref{thm-2-fre-tech}}
Denote $\rho_f:=\rho(\alpha, R_{\varrho_0}e^F)$. Since $|\rho_f-\varrho_0|\leq 2\|F\|_{0}$ by Lemma \ref{Lem-rot-num}, the cocycle $(\alpha, R_{\varrho_0}e^F)$ can be rewritten as $(\alpha, R_{\rho_f}e^{\tilde F})$ for some $\tilde F\in C_{\mathfrak r}^\omega(\T^2, sl(2,\R))$ with $\|\tilde F\|_{\mathfrak r}\leq 30\|F\|_{\mathfrak r}$.
Let $\bar {\mathfrak r}_0:=\mathfrak r, g_0(\phi):=0, F_0(\phi):= \tilde F(\phi), \tilde\varepsilon_0=0$, and $\varepsilon_0=\| F_0\|_{{\mathfrak r}_0}$ satisfy (\ref{condition-varepsilon}). Define ${\mathfrak r}_n, \bar {\mathfrak r}_n, \varepsilon_n, \tilde\varepsilon_n$ as in (\ref{para-induc-2-fre}). By the assumption, cocycle $(\alpha, R_{\rho_f}e^{F_0})\in \tilde {\mathcal{F}}_{{\mathfrak r}_0}(0, \varepsilon_0, \gamma,\tau)$ and applying Lemma \ref{lemma-2-fre-redu} we get the cocycle $(\alpha, R_{\rho_f +\frac{\bar g_0}{2\pi}}e^{\bar F_0})\in \tilde {\mathcal F}_{{\mathfrak r}_1}(\tilde \varepsilon_1, \varepsilon_1, \gamma, \tau)$ by conjugation map $\Psi_0$. Let $g_1:=\bar g_0, F_1:=\bar F_0, \Phi_0:=\Psi_0$. Then the iterative lemma Proposition \ref{prop-2-freq} applies.
Inductively, we obtain the sequence of transformations $\{\Phi_n\}_{n\geq 0}$ with $\|\Phi_n - I\|_{{\mathfrak r}_{n+1}}<\varepsilon_n^{1/3}$ that transforms (\ref{sys-2-fre-before}) to (\ref{sys-2-fre-after}).

Let
\begin{eqnarray*}
\Phi^{(n)}=\Phi_{n-1}\cdots \Phi_0,\  \Phi=\lim_{n\rightarrow \infty}\Phi^{(n)}, \ \  g_\infty=\lim_{n\rightarrow\infty}g_n(\phi).
\end{eqnarray*}
Then under the transformation $x_\infty=\Phi x$, we get the system $(\alpha, R_{\varrho+\frac{g_\infty}{2\pi} })$.  Moreover, we can verify that $\Phi$ is $C^\infty$ in $\phi$: First, by standard computation we get the estimate
$$\|\Phi^{(n+1)}-\Phi^{(n)}\|_{{\mathfrak r}_{n+1}}<2\varepsilon_n^{1/3}.$$
Furthermore, by the selection of $(\alpha, \tilde \alpha)$, for any $j\in\Z^2_+$ there exists $n_j\in\N$ such that for any $n\geq n_j$, we have
$\left(\frac{32}{{\mathfrak r}_0}\tilde q_{n_*+n}^{2\chi}\right)^{|j|}<e^{\frac{\tilde q_{n_*+n}^{1/10}{\mathfrak r}_0}{16}}<\varepsilon_n^{-1/8}.$
Then for all $j\in\Z_+^2, n\geq n_j$, by Cauchy estimates we get
\[\left|D^j(\Phi^{(n+1)}-\Phi^{(n)})\right|\leq \frac{\|\Phi^{(n+1)}-\Phi^{(n)}\|_{{\mathfrak r}_{n+1}}}{{\mathfrak r}_{n+1}^{|j|}}<2\varepsilon_n^{1/3}\cdot \left(\frac{32}{{\mathfrak r}_0}\tilde q_{n_*+n}^{2\chi}\right)^{|j|}<\varepsilon_n^{1/6}. \]
Therefore, $\Phi$ is $C^\infty$ in $\phi$ and hence $g_\infty\in C^\infty(\T^2,\R)$.

\qed\medskip

\section{Appendix: Proof of Proposition \ref{homo-lemma} }

We  solve the approximate equation of (\ref{homo-equ-f-pre})
\begin{equation}\label{homo-equ-approx}
e^{4\pi i\varrho}h(\cdot+\alpha)-h+\mathcal{T}_K(\widetilde{g}(\varphi)h(\cdot+\alpha))=-\mathcal{T}_Kf.
\end{equation}
Let $$ f(\varphi,\theta)=\sum_{l}f_{l}(\varphi)e^{2\pi i\la l,\theta\ra},\qquad
f_{l}(\varphi)=\sum_{k}\hat f_{l}(k)e^{2\pi i k\varphi},$$
\begin{displaymath}
h(\varphi,\theta)=\sum_{|l|<K}h_{l}(\varphi)e^{2\pi i\la l,\theta\ra},\qquad
h_{l}(\varphi)=\sum_{|k|<K-|l|}\hat h_{l}(k)e^{2\pi i
k\varphi}.
\end{displaymath}
The approximate equation (\ref{homo-equ-approx}) is equivalent to (for $|l|<K$)
\begin{eqnarray}
\label{mat}
&&e^{2\pi i (\la l,\alpha' \ra +2\varrho)} h_l(\varphi+\tilde{\alpha})-h_l(\varphi)\\
&&+ e^{2\pi i \la l,\alpha' \ra } \mathcal{T}_{K-|l|}\left(\widetilde{g}(\varphi)h_{l}(\varphi+\tilde{\alpha})\right)
 =-\mathcal{T}_{K-|l|}f_{l}(\varphi).\nonumber
\end{eqnarray}
For any fixed $l$, equation $(\ref{mat})$ can be viewed as a matrix equation
$$(D_{l}+G_{l})\overline{h}_{l} = -\overline{f}_{l},$$
where
$$\overline{h}_{l}=(\hat h_{l}(k))^{T}_{|k|< K-|l|}, \qquad  \overline{f}_{l}=(\hat f_{l}(k))^{T}_{|k|< K-|l|}, $$
and
$$D_l=diag\left(e^{2\pi i (\la l,\alpha' \ra +2\varrho+ k\tilde\alpha)}-1\ :\ |k|< K-|l|\right), $$
$$  G_{l} = \left( e^{2\pi i (\la l,\alpha' \ra+ q\tilde\alpha ) }\hat{\widetilde{g}}(p-q)\right)_{|p|,|q|< K-|l|}.$$

If we denote
$$\Omega_{l,r'}=diag(\cdots, e^{2\pi|k|r'}, \cdots)_{|k|<
K-|l|}$$
 for any $r'\leq r$, then
$$\Omega_{l,r'}({D}_{l}+{G}_{l})\Omega_{l,r'}^{-1}\Omega_{l,r'}\overline{h}_{l} =  -\Omega_{l,r'}\overline{f}_l.$$
Rewrite it as
\begin{eqnarray}\label{h-f-equation}
({D}_{l}+\widetilde{G}_{l,r'})\widetilde{h}_{l,r'} =
-\widetilde{f}_{l,r'},
\end{eqnarray}
where
$$\widetilde{G}_{l,r'}=\Omega_{l,r'}{G}_{l}\Omega_{l,r'}^{-1},\qquad \widetilde{h}_{l,r'}=\Omega_{l,r'}\overline{h}_{l}, \qquad \widetilde{f}_{l,r'}=\Omega_{l,r'}\overline{f}_l, $$ with the estimation
\begin{eqnarray}
\|\widetilde{G}_{l,r'}\|\leq \frac{4\|\widetilde{g}\|_{r}}{1-e^{-2\pi (r-r')}}\leq \frac{16\pi \|g\|_{r}}{1-e^{-2\pi (r-r')}},
\end{eqnarray}
where the matrix norm is the 1-norm of a matrix.
Moreover, one can check easily that
\begin{eqnarray}\label{h-h-tilde}
\|\widetilde{h}_{l,r'} \|=\sum_{|k|<K-|l|}|\hat
h_l(k)|e^{2\pi|k|r'}.
\end{eqnarray}

As $\varrho\in D_{\alpha}(\gamma',\tau')$, then for any $|k|+|l|< K$  we have
$\|D^{-1}_{l}\|\leq\frac{K^{\tau'}}{\gamma'}$
and hence
\begin{eqnarray}
\|D_{l}^{-1}\widetilde G_{l,r-\sigma}\|\leq\frac{16 \pi K^{\tau'}\eta}{\gamma'\sigma}<\frac{1}{2}
\end{eqnarray}
 by (\ref{dia-2}) for $\sigma\leq\frac{1}{4}$.
Thus, there is
$$\|(I+ D_{l}^{-1}\widetilde G_{l,r-\sigma})^{-1}\|<2.$$
Then by (\ref{h-f-equation}), $\widetilde{h}_{l,r-\sigma}$  can be expressed as:
\[\widetilde{h}_{l,r-\sigma}
=-(I+ D_{l}^{-1}\widetilde G_{l,r-\sigma})^{-1}D_{l}^{-1} \widetilde{f}_{l,r-\sigma}.\]
Let $\delta\geq \sigma$. We furthermore use (\ref{h-h-tilde}) to obtain the  estimation
\begin{eqnarray*}
\lefteqn{\|h(\varphi,\theta)\|_{
r-\sigma, s-\delta}}\\&\leq&\sum_{|k|+|l|<K}|\hat
h_l(k)|e^{2\pi|k|(r-\sigma)}e^{2\pi|l|(s-\delta)}\\
&\leq&\sum_{|l|<K}\|\widetilde
h_{l,r-\sigma}\| e^{2\pi|l|(s-\delta)}\\
&\leq&\sum_{|l|<K}\|(I+ D_{l}^{-1}\widetilde
G_{l,r-\sigma})^{-1}\|\cdot\| D_{l}^{-1}\|\cdot\|\widetilde
f_{l,r-\sigma}\| e^{2\pi|l|(s-\delta)}\\
&\leq&\frac{2 K^{\tau'}}{\gamma'}\sum_{|l|<K}\sum_{|k|<K-|l|}\|f\|_{r,s}e^{-2\pi|k|\sigma}e^{-2\pi|l|\delta}\\
&\leq&\frac{C_0(d)K^{\tau'}\widetilde \eta}{\gamma' \sigma^{d}}.
\end{eqnarray*}
Moreover, we can get the control of the error term:
\begin{eqnarray*}
\|\tilde{P}\|_{ r-2\sigma, s-2\delta}&\leq
&\|\mathcal{R}_{K}f\|_{r-2\sigma, s-2\delta}+
\| \mathcal{R}_K(\widetilde{g}(\varphi)h(\cdot+\alpha))\|_{r-2\sigma, s-2\delta}\\
&\leq& C(d)
K^{d}e^{-2\pi K\sigma}\left(\|f\|_{r,s}+\|\widetilde{g}\|_{r}\cdot\|h\|_{r-\sigma,s-\delta}\right)\\
&\leq& C_1(d)K^{d}\tilde
\eta^2\left(1+\frac{K^{\tau'}}{\gamma'\sigma^{d}}\eta\right).
\end{eqnarray*}

\section*{Acknowledgements}
QZ wants to thank Artur Avila for useful discussions.  JW was supported by ``the Fundamental Research Funds for the Central Universities", No. 30918011336, NNSF of China (Nos. 11601230, 11971233, 11701285), and the Natural Science Foundation of Jiangsu Province, China (No. BK20160816). XH was partially  supported by NNSF of China (Grant
 11371019, 11671395) and Self-Determined Research Funds of Central China Normal University (CCNU19QN078). QZ was support by NNSF grant (11671192, 11771077) and Nankai Zhide Foundation.

\end{document}